# Boundary Behavior of Solutions of a Class of Genuinely Nonlinear Hyperbolic Systems

by

Julian Gevirtz

ABSTRACT. We study the set of boundary singularities of arbitrary classical solutions of genuinely nonlinear $2 \times 2$ planar hyperbolic systems of the form $D_k R_k = 0$, where $D_k$ denotes differentiation in the direction $e^{i\theta_k(R_1, R_2)}$, $k = 1, 2$, and where the defining functions $\theta_k$ satisfy

(i) $\theta_2 = \theta_1 + \frac{\pi}{2}$, and

(ii) $A < |\frac{\partial \theta_k(R)}{\partial R_k}| < B$, for all $R \in \mathbb{R}^2$, $k = 1, 2$

for some positive constants $A$, $B$. We show that for any system of this kind there is a $\tau < 1$ such that for any locally Lipschitz solution $R$ in a smoothly bounded domain $G$, the set of points of $\partial G$ at which $R$ fails to have a nontangential limit has Hausdorff dimension at most $\tau$, and, on the other hand, for any such system for which the $\theta_k \in C^\infty(\mathbb{R}^2)$, we construct a $C^\infty$ solution $R$ on a half-plane $\mathbb{H}$ for which the set of points of $\partial \mathbb{H}$ at which $R$ fails to have a nontangential limit has positive Hausdorff dimension. These results are immediately applicable to constant principal strain mappings, which are defined in terms of a system of this kind for which $\theta_1$ is a linear function of $R_1$ and $R_2$.

## 1. Introduction

For hyperbolic systems in two independent variables $x$ and $t$, most often associated with space and time, one usually studies the Cauchy problem in which one seeks a solution $u(x, t)$, $t \geq 0$ for which $u(x, 0)$ coincides with a given $u_0(x)$, the questions considered including well-posedness, global existence, blow-up and behavior of solutions as $t \to \infty$. In the nonlinear case discussion is often limited to initial data with a small range and even for such data, generalized solutions must be considered.

In this paper we concern ourselves with the following inverse question for a certain family of genuinely nonlinear $2 \times 2$ hyperbolic systems: What can be said about the boundary values of an *arbitrary* classical solution in a plane domain $G$? Here "classical" can be taken to mean $C^\infty$, although the treatment we give will be valid for locally Lipschitz solutions. In the first place, we are interested in systems whose formulation imposes no a priori limit on the range of characteristic directions, that is, systems such that for a characteristic given parametrically by $z(s)$, $\arg\{z'(s)\}$ can potentially cover all of $\mathbb{R}$, in contrast to what is implicitly the case in the standard space-time context. Secondly, we are interested in statements valid for *all* solutions rather than ones known to arise from some form of initial value problem. Because of this generality, even in geometrically simple domains such as disks or half-planes characteristics can be quite contorted curves. Although the specific focus of this paper is the size of the set of boundary points at which arbitrary solutions can fail to have nontangential limits, it would be reasonable to investigate other aspects of their behavior and that of the associated characteristics. In any event, given the nonstandard nature of the boundary value question and of several of the issues that arise in dealing with it, we shall begin with a somewhat detailed discussion of a system for which it is physically meaningful, namely the system which describes smooth planar mappings with constant principal stretches (cps-mappings), about which we have previously written ([ChG],[G1]-[G5]). It is in fact the study of the boundary behavior of such mappings that is the main goal of this paper, and we have only chosen to work in a wider context because it is possible to do so with little additional effort, and because this broader approach suggests some interesting questions.

A *cps-mapping* with principal stretch factors $m_1 \neq m_2$ is a mapping $f : G \to \mathbb{C}$ with locally Lipschitz continuous Jacobian $J_f = T(-\phi)\sigma(m_1, m_2)T(\theta)$, where

$$T(\theta) = \begin{bmatrix} \cos\theta & \sin\theta \\ -\sin\theta & \cos\theta \end{bmatrix} \text{ and } \sigma(m_1, m_2) = \begin{bmatrix} m_1 & 0 \\ 0 & m_2 \end{bmatrix}.$$

As is explained in the cited references, apart from regularity considerations, functions $\theta$ and $\phi$ will correspond to such a mapping on a simply connected domain $G$ if and only if they satisfy the autonomous quasilinear hyperbolic system

(1.1) $$D_1(m_1\theta - m_2\phi) = 0; \quad D_2(m_2\theta - m_1\phi) = 0,$$

where

$$D_1 u = (\cos\theta)u_x + (\sin\theta)u_y \text{ and } D_2 u = (-\sin\theta)u_x + (\cos\theta)u_y.$$

The characteristics of a solution are the integral curves of the fields $e^{i\theta}$ and $ie^{i\theta}$. It turns out that a net $\mathcal{N}$ made up of two mutually orthogonal families of curves covering a simply connected $G$ is the net of characteristics of a cps-mapping if and only if for any two curves $C_1, C_2$ belonging to one of the families of $\mathcal{N}$, the change in the inclination of the tangent is the same along all subarcs of curves of the *other* family which join $C_1$ to $C_2$. Nets with this property are known as *Hencky-Prandtl nets*. (See [CS], [G3], [Hem], [Hen], [Hi], [Pr].) The theory of cps-mappings we have developed is based on direct application of (1.1) together with this Hencky-Prandtl (HP) property and the equations

(1.2) $$D_2 D_1 \theta = [D_1\theta]^2 \quad \text{and} \quad D_1 D_2 \theta = -[D_2\theta]^2,$$

which are also effectively equivalent to (1.1) and which are very special cases of equations derived by Lax [L] in the context of considerably more general genuinely nonlinear $2 \times 2$ hyperbolic systems in the plane and used by him in connection with the inevitability of singularity formation. The *blow-up equations* (1.2) imply that if a characteristic $C$ has curvature $\kappa_0$ at $p$ then the orthogonal characteristic arc emanating from $p$ towards the concave side of $C$ can have length at most $\frac{1}{\kappa_0}$, that is, that the boundary of $G$ must be encountered after moving at most a distance of $\frac{1}{\kappa_0}$ along this orthogonal characteristic. A *characteristic length bound* of this kind is a *sine qua non* for the theory we are developing and plays a fundamental role in all that is to follow.

When regarded as deformations with constant principal strains, cps-mappings are of concrete interest as models in a number of physically interesting contexts (see [Y]). Consider, for example, a thin liquid film on a plane surface which upon solidification takes on a rectangular cryptocrystalline structure, that is, at each point a suitably oriented minute square of the original liquid becomes a rectangular crystal whose side lengths are constant multiples of the side length of the square. In this light global geometric results for cps-mappings tell one about the extent to which the shape of the original film can change as a result of such solidification and about how matter is moved around in the process, and statements about the existence of boundary limits of $\theta$ (and, in light of (1.1), of $\phi$, and consequently of the Jacobian of the mapping) tell one to what extent the cryptocrystalline structure is present at the very edge of the solidified lamina. Applied to the system (1.1) the main result of this paper says that *there is some number $\tau < 1$ such that if $G$ is smoothly bounded, then $\theta$ and $\phi$ can fail to have nontangential limits on a set $S \subset \partial G$ of Hausdorff dimension at most $\tau$*. On the other hand, the construction of Section 5 shows that *the set of boundary points at which $\theta$ does not have nontangential limits can in fact have positive Hausdorff dimension.*

Beyond their immediate physical significance, cps-mappings constitute a particularly important and tractable class of planar quasi-isometries, for which we believe they will ultimately be shown to display extremal behavior for many of the as yet unsolved distortion questions (see [J1], [J2]). In this direction a very significant inroad was made a few years ago by Gutlyanskii and Martio [GM], who showed that the *spiral mappings* given by

$$G(re^{i\psi}) = re^{i(\psi + \psi_0 \frac{\log r}{\log \rho})},$$

are extremal for the problem of determining for given $\rho > 1$ and $\psi_0 > 0$, the smallest ratio $\frac{m_1}{m_2} > 1$, such that there is a quasi-isometry $f$ of the annulus $1 < |z| < \rho$ onto itself with stretching bounds $m_1$, $m_2$ which satisfies the boundary conditions

$$f(z) = z \quad \text{and} \quad f(\rho z) = \rho e^{i\psi_0} z, \quad \text{for } |z| = 1.$$

It is in fact not hard to see that $G$ is indeed a cps-mapping of the entire punctured plane $\mathbb{C}\setminus\{0\}$ with

$$m_1 = \tfrac{\sqrt{a^2+4}+a}{2} \quad \text{and} \quad m_2 = \tfrac{\sqrt{a^2+4}-a}{2},$$

where $a = \frac{\psi_0}{\log \rho}$. We call these spiral mappings because the corresponding inclination functions are of the form

(1.3) $$\theta(re^{i\psi}) = \psi + \alpha,$$

where $\tan \alpha = m_2$, which is to say that the characteristics form two mutually orthogonal families of logarithmic spirals all members of each of which are rotations of each other.

We now describe the class of $2 \times 2$ systems for which we treat the boundary singularity question. Let $\theta_k = \theta_k(R_1, R_2), k = 1, 2$. For given $R_1(x,y)$, $R_2(x,y)$ and any $u = u(x,y)$ we write

(1.4) $\quad D_k u = \cos\theta_k(R_1(x,y), R_2(x,y)) \frac{\partial u}{\partial x} + \sin\theta_k(R_1(x,y), R_2(x,y)) \frac{\partial u}{\partial y}.$

It is well-known that in general an autonomous $2 \times 2$ quasi-linear homogeneous hyperbolic system for unknown functions $f$ and $g$ is formally equivalent to a system of the form

$$D_k R_k = 0, \ k = 1, 2,$$

with appropriate inclination functions $\theta_k(R_1, R_2)$. The relationship between the *Riemann invariants* $R_1$, $R_2$ and $f$, $g$ is of the form $R_k(x,y) = F_k(f(x,y), g(x,y))$. Henceforth we write $\Theta = (\theta_1, \theta_2)$ and use the term *system* to mean a $C^\infty$ mapping $\Theta : \mathbb{R}^2 \to \mathbb{R}^2$, although the smoothness requirement could be weakened substantially. Obviously, a *solution* of the system $\Theta$ in a domain $G$ of the plane is then a pair of functions $R_1(x,y)$, $R_2(x,y)$ for which $R_k$ is constant on each integral curve (henceforth referred to as a *k-characteristic*) of the field $e^{i\theta_k(R_1(x,y), R_2(x,y))}$, $k = 1, 2$. A system is said to be *genuinely nonlinear* if the derivatives $\frac{\partial \theta_k}{\partial R_k}$, $k = 1, 2$, never vanish. It is clear that the system (1.1) for the $\theta$ and $\phi$ associated with cps-mappings is already in Riemann invariant form with $R_i = m_i \theta - m_j \phi$, so that in this case two inclination functions are given by

(1.5) $\quad \theta_1 = \theta = \frac{m_1 R_1 - m_2 R_2}{m_1^2 - m_2^2} \quad \text{and} \quad \theta_2 = \theta + \frac{\pi}{2}.$

(Here we have used the convention, in force throughout this paper, to the effect that $\{i, j\} = \{1, 2\}$.) This system is obviously genuinely nonlinear and in fact is the simplest possible such system in that the two families of characteristics are mutually orthogonal and $\Theta$ is a linear function of $R = (R_1, R_2)$. We now define the family of systems with which we work.

**Definition 1.1**. A *normal system* is a system $\Theta$ for which
(*i*) $\theta_2 = \theta_1 + \frac{\pi}{2}$ and
(*ii*) There are constants $A$, $B > 0$ such that $A < |\frac{\partial \theta_k}{\partial R_k}| < B$ for all $R \in \mathbb{R}^2$.

The only hyperbolic systems with which we deal will be normal systems, for which we use the symbol $\theta$ to denote $\theta_1$. We shall show that *for any normal system $\Theta$ there is a $\tau = \tau(\Theta) < 1$ such that for any smoothly bounded $G \subset \mathbb{C}$ and any solution $R$ of $\Theta$ on $G$ the set of points of $\partial G$ at which $R$ does not have a nontangential limit has Hausdorff dimension at most $\tau$.* This will follow as an immediate consequence (Corollary 4.2) of our principal result, Main Theorem 2.2, which deals with boundary singularities of a class of functions effectively more general than the class of solutions of normal systems. An outline of the proof, which is actually made simpler by this somewhat greater generality, is given early in Section 2, just after the statement of the main theorem. Furthermore, in Section 5 we shall show that *for any such $\Theta$ there is a $C^\infty$ solution $R$ in the upper half-plane $\mathbb{H}$ for which the set of points of $\partial \mathbb{H}$ at which $R$ does not have nontangential limits has positive Hausdorff dimension.*

## 2. Quasi-HP Functions, the Characteristic Length Bound and Related Matters

Let $G \subset \mathbb{C}$ be a domain. If $\theta$ is a locally Lipschitz function on $G$, the integral curves of the fields $e^{i\theta(z)}$ and $ie^{i\theta(z)}$ will be called the 1- and 2-characteristics of $\theta$, respectively. The term *full characteristic* will refer to the complete integral curve and we will use the term *half-characteristic* to refer to either of the two arcs into which a full characteristic is divided by one of its points. (Characteristics which are closed curves will not arise in this paper.) As indicated in the Introduction we will use the convention $\{i, j\} = \{1, 2\}$ throughout. Arcs of $k$-characteristics will be called $k$-arcs, or less specifically, characteristic arcs. With reference to a given $\theta$, a characteristic arc joining points $a, b \in D$ will be denoted by $ab$ and we shall use the abbreviation
$$\Delta\theta(ab) = \theta(b) - \theta(a).$$
A domain $Q \subset G$ will be said to be a *positively* (*negatively*) *oriented characteristic quadrilateral* of $\theta$, and we write $Q = abcd$, if $\partial Q$ is a Jordan curve lying in $G$ containing four points $a$, $b$, $d$, $c$ occurring in that order when $\partial Q$ is traversed in the positive (negative) sense and such that $ab$ and $cd$ are $i$-arcs and $ac$ and $bd$ are $j$-arcs. We say that $ab$ and $cd$ are *translates* of each other with respect to or along any $j$-characteristic passing through $ab$. For a curve parametrized by $z = z(s)$, $\alpha < s < \beta$, we use the terms "to the right of $C$", "to the left of $C$" in the obvious sense, so that, for example, if $C$ is a characteristic arc and $p \in C$, we describe an orthogonal characteristic arc or a half-characteristic as emanating from $p$ to the right or left of $C$.

We shall denote the 2-dimensional measure of $X \subset \mathbb{C}$ by $\mu(X)$ and the 1-dimensional measure of a set $A$ by $\lambda(A)$. The parameter $s$ will always refer to arc length. We use the notation $N(a, r) = \{z : |z - a| < r\}$ and denote the line segment joining $a$ to $b$ by $\overline{ab}$. The overline will also be used to denote closure, but this should cause no confusion. For $X, Y \subset \mathbb{C}$, $\text{dist}(X, Y) = \inf\{|y - x| : x \in X, y \in Y\}$, and for $a \in \mathbb{C}$, $\text{dist}(a, X) = \text{dist}(\{a\}, X)$.

**Definition 2.1.** Let $K \geq 1$. A locally Lipschitz function $\theta$ on a domain $G$ will be said to be a $K$-*quasi-HP function* or to have the $K$-*quasi-HP property* if for any characteristic quadrilateral $abcd \subset G$ there holds
(2.1) $$\tfrac{1}{K}|\Delta\theta(ac)| \leq |\Delta\theta(bd)| \leq K|\Delta\theta(ac)|.$$
The net consisting of all of the full 1- and 2-characteristics of a $K$-quasi-HP function will be called a $K$-*quasi-HP net*.

A simple continuity argument shows that this definition implies that the $\Delta\theta(ac)$ and $\Delta\theta(bd)$ in (2.1) must in fact have the same sign (unless both vanish). It is also obvious that a 1-quasi-HP net is an HP-net. Note that while the HP-property is a local condition that implies its global counterpart, this is not the case for the $K$-quasi-HP property when $K > 1$. *The net of characteristics of any locally Lipschitz solution to a normal system $\Theta$ is a $K$-quasi-HP net*, where $K = K(\Theta)$. Indeed, if $A$ and $B$ are as in Definition 1.1, then it is clear that we can take $K(\Theta) = B/A$.

We can now state our

**Main Theorem 2.2**. *There is a number $\tau = \tau(K) < 1$ with the property that if $\theta$ is any $K$-quasi-HP function on a Jordan domain $G$ with $C^2$ boundary, then the set of points of $\partial G$ at which $\theta$ fails to have a nontangential limit has Hausdorff dimension at most $\tau$.*

Because the proof of this theorem, to be given in Section 4, is quite involved and depends on the prior development of a considerable amount of machinery in this and the following section, we shall briefly explain here how it proceeds. As we shall show (see Corollary 3.29) for any point $p \in \partial G$ at which $\theta$ does not have a nontangential limit there is either a nontrivial fan of characteristics emanating from $p$, or for $k = 1$ or 2, $p$ is not the endpoint of a $k$-characteristic; in the latter case we say that $p$ is a $k$-singularity. Since, as will be apparent, any quasi-HP function can have at most a countable number of fans, we only need to show that the Hausdorff dimension, $\dim(S)$, of the set $S$ of $k$-singularities satisfies $\dim(S) \leq \tau < 1$. For this to be the case it is enough that there be some $\delta = \delta(K) > 0$ such that any almost straight arc $A \subset \partial G$ has a subarc of length at least $\delta\lambda(A)$ which has at most a countable set of $k$-singularities of $\theta$. If for some $\theta$ there were no such $\delta$, then for any $N$ there would have to be an almost straight arc $A \subset \partial G$ such that there is a set $S_N$ of $k$-singularities of $\theta$ which are essentially uniformly distributed along $A$. It is not too hard to show that any $k$-singularity $p$ is surrounded by arbitrary small $k$-characteristics which join boundary points on either side of $p$ (see Proposition 3.19). For sufficiently large $N$, starting with a set of $N$ small $k$-characteristics, one surrounding each of the points of $S_N$, we show that there must be a $k$-characteristic $C$ (where the term "$k$-characteristic" is used in an appropriate sense - see the discussion of *extended characteristics* between Propositions 3.13 and 3.14) whose endpoints lie on $A$ and are at least $\delta'\lambda(A)$ apart (where $\delta' > 0$ depends solely on $K$) and which is tangent to $A$ at a point $m \in A$, where $m$ is appropriately bounded away from the endpoints of $C$. We then use this to obtain a subarc $A'$ of $A$ which contains $m$, whose length is bounded below by $\delta''\lambda(A)$, where $\delta''$, like $\delta'$, depends only on $K$, and on which $\theta$ can have at most a countable number of $k$-singularities (see Proposition 3.33), thereby arriving at the desired contradiction. A good measure of the complexity of the proof lies in establishing the existence of the extended characteristic $C$, which is carried out in Section 4, but which depends on properties of the net of extended characteristics of quasi-HP-functions developed in Section 3.

We begin with the following proposition, which is immediate.

**Proposition 2.3.** (Invariance of $K$) *Let $\theta$ be a $K$-quasi-HP function on $G$, and $a, b \in \mathbb{C}$ with $a \neq 0$. Then $\theta(az + b) - \arg a$ is a $K$-quasi-HP function on $\frac{1}{a}(G - b)$.*

The next proposition is a special case of [G1, Lemma 2]; its proof is included for the sake of completeness. A function $\theta$ is said to be a locally $L$-Lipschitz function on $G$ if each point of $G$ has a neighborhood in which $\theta$ satisfies a Lipschitz condition with constant $L$.

**Proposition 2.4.** (Length Change Estimate) *Let $G$ be a simply connected domain and let $\theta$ be a locally $L$-Lipschitz function on $G$. Let $Q = abcd \subset G$ be a positively oriented characteristic quadrilateral such that $|b - a| = l$, $|c - a| = \epsilon \leq l$, and $\mathrm{dist}(Q, \partial G) = \eta > 0$. There is some $\bar{l} = \bar{l}(L, \eta) > 0$ such that*
$$(2.2) \qquad |d - c| = l - \epsilon \Delta\theta(ab) + O(\epsilon^2 + l^3) \quad \text{for } l \leq \bar{l},$$
*where the constant implied by the big-O depends only on $L$ and $\eta$.*

**Proof.** In what follows, when we say that some quantity is big-$O$ of some expression, we mean that this is so for all $l$ less than some positive number which depends only on $L$ and $\eta$, and that the constant corresponding to the big-$O$ depends only on $L$ and $\eta$. Let $ab$ and $cd$ be $i$-arcs. Without loss of generality we can assume that $a = 0$ and $b = l$. Let $\frac{\pi}{2} + \alpha$ and $\frac{\pi}{2} + \beta$ be the inclinations of the tangents to the $j$-characteristics at $a$ and $b$, respectively, and let $\frac{\pi}{2} + \alpha' = \arg\{c - a\}$, $\frac{\pi}{2} + \beta' = \arg\{d - b\}$. Clearly, $\alpha, \beta, \alpha'$ and $\beta'$ are all $O(l)$ and $\beta - \alpha = \theta_0$. It easily follows from the Lipschitz condition that $\alpha' = \alpha + O(\epsilon)$. We have $d = \epsilon i e^{i\alpha'} + t e^{i\gamma}$, where $t = |d - c| = O(l)$ and $\gamma = O(l)$. We also have $d = l + s i e^{i\beta'}$, where $s = O(l)$. Just as $\alpha' = \alpha + O(\epsilon)$, one sees that $\beta' = \beta + O(s)$. From the two expressions for $d$ we have that
$$\epsilon i e^{i\alpha'} + t e^{i\gamma} = l + s i e^{i\beta'},$$
so that considering real and imaginary parts we have
$$(2.3) \qquad -\epsilon \sin \alpha' + t \cos \gamma = l - s \sin \beta'$$
and
$$\epsilon \cos \alpha' + t \sin \gamma = s \cos \beta'.$$
The latter equation implies that
$$s(1 + O(l^2)) = \epsilon(1 + O(l^2)) + O(l^2),$$
so that $s = \epsilon + O(l^2)$. Thus, since $\beta' = \beta + O(s)$, we have
$$(2.4) \qquad \beta' = \beta + O(\epsilon + l^2).$$
From (2.3) it now follows that
$$t(1 + O(l^2)) = l + \epsilon \sin \alpha' - s \sin \beta' = l + \epsilon(\alpha' + O(l^3)) - (\epsilon + O(l^2))(\beta' + O(l^3)),$$
so that from the fact that $\alpha' = \alpha + O(\epsilon)$ and (2.4) it follows that
$$|d - c| = t = l + \epsilon(\alpha' - \beta') + O(l^3)$$
$$= l + \epsilon(\alpha + O(\epsilon)) - \epsilon(\beta + O(\epsilon + l^2)) + O(l^3)$$
$$= l - \epsilon(\beta - \alpha) + O(\epsilon^2 + l^3).$$
Since $\beta - \alpha = \Delta\theta(ab)$, we are done. ∎

Let $C$ be a characteristic arc parametrized by $z = z(s)$, $\alpha < s < \beta$. Let $D^+(\theta, C)$ denote the infimum of all $\tau > 0$ such that $\frac{\theta(s_1) - \theta(s_2)}{s_1 - s_2} \leq \tau$ for all $s_1, s_2 \in (\alpha, \beta)$. Note that $D^+(\theta, C)$ depends on the orientation of $C$ implicit in the parametrization.

**Proposition 2.5.** (Characteristic Length Bound). *Let $\theta$ be a $K$-quasi-HP function on $G$ and let $C \subset G$ be an open $i$-arc with $D^+(\theta, C) = \kappa > 0$. Then for any $\eta > 0$ there is* *some point $p \in C$ such that the $j$-half-characteristic emanating from $p$ to the left of $C$ has length at most $\frac{K}{\kappa} + \eta$.*

**Proof.** Let $\delta < 1$ be a small positive number. Let $C$ be parametrized by $z = z(s)$. Clearly there is a point $p = z(s_0)$ such that $\theta$ is differentiable at $p$ and $\theta'(s_0) \geq \kappa(1 - \frac{\delta}{2})$. Without loss of generality we can take $s_0 = 0$. Let $J$ be the $j$-half-characteristic emanating from $p$ to the left of $C$. Let $C'$ be any compact arc of $J$ whose initial point is $p$. It is clearly sufficient to show that

(2.5) $$\lambda(C') \leq \frac{K}{\kappa(1-\delta)^3} + \delta.$$

Let $t_0 < \delta$ be so small that the entire characteristic quadrilateral $Q$ for which an $i$-side is $z([0, t_0])$ and a $j$-side is $C'$ is contained in $G$. Let $2\eta = \text{dist}(Q, \partial G) > 0$. Let $L$ be such that $\theta$ is an $L$-Lipschitz function in the $\eta$-neighborhood of $Q$. By decreasing the size of $t_0$ if necessary we can assume that

$$\frac{\theta(z(t)) - \theta(p)}{|z(t) - p|} \geq \kappa(1 - \delta), \text{ for } 0 < t \leq t_0$$

and that $t_0 \leq \bar{l}(L, \eta)$, where $\bar{l}(L, \eta)$ is as in the preceding proposition. Let $\theta_1(t) = \theta(z(t)) - \theta(p)$. Clearly, by taking a smaller positive value for $t_0$ we can assume that for all $t \in (0, t_0]$ and all $\theta_0 \geq \frac{\theta_1(t)}{K}$, if $l \leq |z(t) - p|$ and $\epsilon = l^{3/2}$, then

(2.6) $$l - \theta_0 \epsilon + O(\epsilon^2 + l^3) \leq l - \theta_0(1 - \delta)\epsilon,$$

where the big-$O$ is that of the previous proposition. Let $C'$ be parametrized by $w = w(s)$, $0 \leq s \leq \lambda_0$. We take $t \leq t_0$ such that if $a_0, a_1, \cdots, a_n$ is an ordered sequence of points on $C'$, the distance between any successive two of which is at most $t$, then

(2.7) $$\sum_{k=1}^{n} |a_k - a_{k-1}| \geq (1 - \delta)\lambda(a_0 a_n).$$

Let $l_0 = |z(t) - p|$, and let $a_0 = p$, $b_0 = z(t)$. Let $C''$ be the other $j$-side of $Q$. Let $d = \text{dist}(C', C'') > 0$. Let $a_1 \in C'$ to be the (first) point for which $|a_1 - a_0| = \epsilon_0 = l_0^{3/2}$. Let $C_1$ be the $i$-half-characteristic emanating from $a_1$ to the right of $J$ and let $b_1$ be the point common to $C_1$ and $C''$. From Proposition 2.4 and (2.6) it follows that

$$l_1 = |b_1 - a_1| \leq l_0 - \theta_1(t)(1 - \delta)\epsilon_0.$$

Clearly, we can continue this process until we come to an $a_n \in C'$ within $l_0$ of $w(\lambda_0)$. Taking into account that $\theta$ is a $K$-quasi-HP function and (2.7) we see that for this $n$ we have that

$$0 < d \leq |b_n - a_n| \leq l_0 - \frac{\theta_1(t)}{K}(1-\delta)\sum_{k=1}^{n}|a_k - a_{k-1}| \leq l_0 - \frac{\theta_1(t)}{K}(1-\delta)^2 \lambda(a_0 a_n)$$

$$\leq l_0 - \frac{\theta_1(t)}{K}(1-\delta)^2(\lambda(C') - l_0) \leq l_0 - \frac{\theta_1(t)}{K}(1-\delta)^2(\lambda(C') - \delta),$$

so that $\lambda(C') - \delta < \frac{l_0}{\frac{\theta_1(t)}{K}(1-\delta)^2} = \frac{K}{\frac{\theta_1(t)}{l_0}(1-\delta)^2} \leq \frac{K}{\kappa(1-\delta)^3}$. Thus indeed (2.5) holds. ∎

If $z = z(s)$, $\alpha < s < \beta$, is a parametrization of an $i$-arc $C$, then $\kappa_i(z(s)) = d\theta(z(s))/ds$ exists almost everywhere on $(\alpha, \beta)$ and gives the curvature of $C$ at $z(s)$. It follows immediately from the preceding proposition that if $C$ is an $i$-arc and $\kappa_i(p)$ exists, then the $j$-half-characteristic emanating from $p$ towards the *concave* side of $C$ (that is, to the left or right of $C$ according as $\kappa_i(p) > 0$ or $\kappa_i(p) < 0$) has length at most $\frac{K}{\kappa_i(p)}$.

Let $\phi$ be a real-valued function defined on $(\alpha, \beta)$, and let $s_0 \in (\alpha, \beta)$. We denote by $D^-\phi(s_0)$ and $D^+\phi(s_0)$ the lower and upper limits of $\frac{\phi(s) - \phi(s_0)}{s - s_0}$ as $s \to s_0$ and by $D_R^-\phi(s_0)$ and $D_L^-\phi(s_0)$ ($D_R^+\phi(s_0)$ and $D_L^+\phi(s_0)$) the corresponding right and left lower (upper) limits. We have the following simple corollary of the characteristic length bound.

**Proposition 2.6.** (Curvature Bound). *Let $z = z(s)$, $\alpha < s < \beta$ be an $i$-characteristic of a $K$-quasi-HP function $\theta$ on $G$ and let $J$ be the $j$-half-characteristic emanating from $z(s_0)$ to the left of $C$ and joining it to $\partial G$. Then $D^+\theta(z(s_0)) \leq \frac{K}{\lambda(J)}$.*

**Proposition 2.7.** (Length Monotonicity). *Let $abcd$ be a positively oriented characteristic quadrilateral of a quasi-HP function $\theta$ on $G$. Let $z = z(s)$, $\alpha \leq s \leq \beta$, $z(\alpha) = a$, parametrize $ab$. Let $E$ be a $\lambda$-measurable subset of $ab$ such that $\theta'(s) \leq 0$ at all points of $E$. (That is, $ab$ is nonconcave towards the inside of $abcd$ at all points of $E$). If $E'$ is the set of points of $cd$ which are joined to points of $E$ by translates of $ac$, then $\lambda(E') \geq \lambda(E)$.*

**Proof**. This is a simple consequence of the length change estimate, the fact that $\theta$ is a quasi-HP function on $G$ and elementary measure theory.∎

In Section 4 we shall make use of the following two lower bounds for the area of certain regions made up of families of characteristic arcs of a $K$-quasi-HP function $\theta$. They follow as simple corollaries of Proposition 2.6; the constants, $\eta$ and $\eta'$ which appear in them depend solely on $K$.

**Proposition 2.8.** (Area Bound 1). *If $C$ is the interior of an $i$-arc of a $K$-quasi-HP function $\theta$ of length $\lambda_i$, and for each $w \in C$, $C'(w)$ is a $j$-arc of $C$ containing $w$ and of length at least $\lambda_j$ such that $U = \cup \{C'(w) : w \in C\}$ is open, then $\mu(U) \geq \eta\lambda_i\lambda_j$.*

**Proposition 2.9.** (Area Bound 2). *If $C$ is the interior of an $i$-arc of a $K$-quasi-HP function $\theta$, and for each $w \in C$, $C'(w)$ is a $j$-arc of $C$ containing $w$ and of length at least $\lambda_j$ such that $U = \cup \{C'(w) : w \in C\}$ is open, then $\mu(U) \geq \eta'\lambda_j^2|\Delta\theta(C)|$.*

### 3. Extended Characteristics, Regular and Singular Boundary Behavior

Our approach to boundary behavior requires the examination of curves which are in effect characteristics whose interiors (i.e, sets of nonendpoints) contain boundary points; the subtleties that arise in this connection require careful discussion. Hereafter the symbol $\mathcal{G}$ will denote the family of all Jordan domains $G \subset \mathbb{C}$ for which $\partial G$ is a $C^2$ curve and $\mathcal{G}(\rho) \subset \mathcal{G}$ will denote the family of all $G \in \mathcal{G}$ such that for each $p \in \partial G$ the interior of one of the circles of radius $\rho$ tangent to $\partial G$ at $p$ is contained in $G$ and that of the other such circle is contained in $\mathbb{C}\backslash\overline{G}$. Obviously, for $G \in \mathcal{G}(\rho)$ the unsigned curvature of $\partial G$ is everywhere bounded by $\frac{1}{\rho}$, so that $\mathcal{G} = \cup \{\mathcal{G}(\rho) : \rho > 0\}$. Furthermore, HP$(G, K)$ will denote the family of $K$-quasi-HP functions on $G$. Although sometimes the hypotheses of the propositions of this section do not state so explicitly, they always deal with $G \in \mathcal{G}$ and $\theta \in \text{HP}(G, K)$. (Although the results of this paper apply to unbounded and multiply connected domains as well, the proof of the main theorem itself will entail only consideration of Jordan domains, and in fact we will be able to work largely with Jordan domains of the kind we call "characteristic subdomains," as defined below.) By an *arc* $C$ we shall henceforth mean a continuous one-to-one mapping $z = z(t)$ of an interval $[\alpha, \beta]$ into the closure $\overline{G}$ of $G$ and $z((\alpha, \beta))$ will be referred to as the *interior* of $C$. As in Section 2, when an arc is considered to be oriented, we use the term "to (towards) the right (left) of $C$" to refer to the part of $C$ (immediately) to the right (left) of $C$, and the term *full characteristic* will refer to a complete integral curve of $\theta$ or $\theta + \frac{\pi}{2}$ in $G$. The following proposition was proved in [G5, Proposition 2.8] for HP-nets. Although the proof is virtually the same in the present more general context we include it here for the sake of completeness.

**Proposition 3.1.** *Let $C$ be a full characteristic of a quasi-HP function in $G$ parametrized by $z = z(s)$, $s \in (\alpha, \beta)$. Then $\lim\limits_{s \to \beta} z(s)$ exists and belongs to $\partial G$.*

**Proof**. Clearly the conclusion holds if $\beta \neq \infty$, so that we assume $\beta = \infty$. First of all, we show that $\text{dist}(z(s), \partial G) \to 0$ as $s \to \infty$. If this were not true, then there would be a $z_0 \in G$ and an $\epsilon > 0$ such that for some sequence $\{s_i\}$ tending to $\infty$, $z(s_i) \to z_0$, but $z([s_i, s_{i+1}]) \cap \partial N(z_0, \epsilon) \neq \emptyset$. But from this it would follow that some orthogonal characteristic crosses $C$ twice, an impossible occurrence in light of the simple connectivity of $G$. We can now show that, in fact, $z(s) \to b \in \partial G$ as $s \to \infty$. If this is not so, the foregoing then implies that there is an arc $E$ of $\partial G$, $\lambda(E) > 0$, each point of which is an accumulation point of $C_\gamma = \{z(s) : s > \gamma\}$ for each $\gamma \in (\alpha, \infty)$. Since $G$ is bounded and $\beta = \infty$, $C$ cannot be a straight line, so that from the characteristic length bound it follows that there is an orthogonal half-characteristic $C'$ of finite length which joins some $z(\sigma)$ to a point $e \in \partial G$. Since $C$ cannot cross $C'$ twice in $G$, $C_\sigma \subset G \setminus C'$. Let $z_1, z_2$ be distinct points of $E \setminus \{e\}$. For each $\delta > 0$, $C_\sigma$ has a subarc $pp' \subset N(\partial G, \delta) \setminus C'$, with $p, p' \in N(z_1, \delta)$ and a point $p'' \in pp' \cap N(z_2, \delta)$. For obvious topological reasons, for each sufficiently small $\delta$, there must be a point $q$ on $pp'$ which is joined to a point in $N(z_1, \delta)$ by an orthogonal characteristic arc $B$ of length at least $|z_1 - z_2| - 2\delta$ such that the curvature of $C$ at $q$ exists and tends to infinity as $\delta \to 0$ and $C$ is concave towards the side from which $B$ emanates. But this clearly violates the characteristic length bound (Proposition 2.5), as indicated in the paragraph immediately following its proof.∎

**Definition 3.2.** An arc $C$ for which $z(\alpha), z(\beta) \in \partial G$ and for which $z((\alpha, \beta)) \subset G$ and is a full $i$-characteristic will be called an *elementary $i$-characteristic*.

In other words, an elementary characteristic is a full characteristic together with its endpoints, which are well defined by the preceding proposition. Note that for each $p \in \partial G$, $\{p\}$ is a *trivial* elementary $i$-characteristic. An $i$-*characteristic arc* (also called an $i$-*arc*) will be a subarc of an elementary $i$-characteristic.

**Lemma 3.3.** *Let $G \in \mathcal{G}$. Then there exists a number $B' = B'(G) \in (0, 1]$ with the following property. Let $a, b \in \partial G$. Let $C_1$ and $C_2$ be the two closed arcs into which $\partial G$ is divided by $a$ and $b$. (If $a = b$, $C_1 = \{a\}$ and $C_2 = \partial G$.) Then*
$$\text{dist}(z, C_1) + \text{dist}(z, C_2) \geq B' \text{dist}(z, \{a, b\}), \text{ for all } z \in G.$$

**Proof**. This is self-evident.∎

**Proposition 3.4.** (Bounded Length of Characteristics). *Let $G \in \mathcal{G}$. There is some $M = M(G, K)$ such that $\lambda(C) \leq M$ for all elementary characteristics $C$ of any $\theta \in \text{HP}(G, K)$.*

**Proof.** Let $\theta$ be a $K$-quasi-HP function on $G$ and let $C$ be an elementary $i$-characteristic of $\theta$. We regard $C$ as being oriented and let $A = C \cap \partial G$ (that is, $A$ is the set of endpoints of $C$). Let $d_0 = \text{diam}(G)$, and let $d = \sup\{\text{dist}(z, A) : z \in C\} \leq d_0$. Clearly,
(3.1) $\qquad \mu(\{z \in G : \text{dist}(z, A) \leq 8r\}) \leq Br^2,$
where $B = 128\pi$. For $k \geq 0$, let
$$G_k = \{z \in G : \tfrac{d}{2^k} \leq \text{dist}(z, A) \leq \tfrac{d}{2^{k-1}}\}$$

and
$$C_k = C \cap G_k,$$
so that $C \cap G = \cup \{C_k : k \geq 1\}$. For each nonendpoint $p$ of $C$, let $J(p)$ denote the elementary $j$-characteristic containing $p$. Obviously, $J(p) \cap J(p') \cap G = \emptyset$ for $p \neq p'$. For $k \geq 1$ and $p \in C_k$ let $l(p)$ and $r(p)$ be the first points encountered when moving along $J(p)$ from $p$ to the left and right of $C$, respectively, which are not in the interior of $G_{k-1} \cup G_k \cup G_{k+1}$. Each $q \in \{l(p), r(p)\}$ is either on $\partial G$ or is in $G$ and satisfies one of $\text{dist}(q, A) = \frac{d}{2^{k+1}}$ or $\text{dist}(q, A) = \frac{d}{2^{k-2}}$. Say $q \notin \partial G$ and $\text{dist}(q, A) = \frac{d}{2^{k+1}}$. Then $|p - q| \geq \frac{d}{2^{k+1}}$, since otherwise $\text{dist}(p, A) \leq |p - q| + \text{dist}(q, A) < \frac{d}{2^k}$, which contradicts the fact that $p \in C_k$. If $q \notin \partial G$ and $\text{dist}(q, A) = \frac{d}{2^{k-2}}$, then $|p - q| \geq \frac{d}{2^{k+1}}$ since otherwise
$$\text{dist}(p, A) \geq \text{dist}(q, A) - |p - q| > \frac{d}{2^{k-2}} - \frac{d}{2^{k+1}} > \frac{d}{2^{k-1}},$$
which is inconsistent with $p \in C_k$. Thus $|p - q| \geq \frac{d}{2^{k+1}}$ if $q \notin \partial G$. Hence, if at least one of $l(p), r(p)$ is not in $\partial G$ we have for the open subarc $J_1(p)$ of $J(p)$ with endpoints $l(p)$, $r(p)$ that $\lambda(J_1(p)) \geq \frac{d}{2^{k+1}}$. If, on the contrary, $\{l(p), r(p)\} \subset \partial G$, then it follows from the preceding lemma that
$$\lambda(J_1(p)) \geq |l(p) - p| + |p - r(p)| \geq \frac{B'd}{2^k} > \frac{B'd}{2^{k+1}},$$
so that this bound holds in all cases for $p \in C_k$ since $B' \leq 1$. Since $J_1(p) \setminus \{l(p), r(p)\} \subset G_{k-1} \cup G_k \cup G_{k+1}$ for $p \in C_k$, it follows from (3.1) that
$$B\left(\frac{d}{2^{k+1}}\right)^2 \geq \mu(\{z \in G : \frac{d}{2^{k+1}} \leq \text{dist}(z, \partial G) \leq \frac{d}{2^{k-2}}\}) = \mu(G_{k-1} \cup G_k \cup G_{k+1})$$
$$\geq \mu(\cup \{J_1(p) : p \in C_k\}) \geq \eta \lambda(C_k) \frac{B'd}{2^{k+1}},$$
by Proposition 2.8. From this we have that $\lambda(C_k) \leq \frac{Bd}{\eta B' 2^{k+1}}$. But since $C = \cup \{C_k : k \geq 1\}$, we conclude that $\lambda(C) \leq \frac{Bd}{2B'\eta} \leq \frac{Bd_0}{2B'\eta}$. Since $B = 128\pi$ and $B'$ and $d_0$ depend only on $G$, and $\eta$ depends only on $K$, we are done. ∎

**Definition 3.5.** An elementary characteristic $C$ one of whose endpoints is $p \in \partial G$ said to *exit* $G$ at $p$. If $C$ is parametrized by $z = z(s)$, $0 \leq s \leq L$ and $\lim_{s \to 0} \theta(z(s))$ exists, then $C$ is said to exit *regularly* at $z(0)$, otherwise it is said to exit *singularly*.

**Proposition 3.6.** *Let $G \in \mathcal{G}$ and $\theta \in \text{HP}(G, K)$. Let $C = ab \subset G$ be a closed $i$-arc of $\theta$ and let $E$ be an arc joining $a$ and $b$ in $G$ such that $E \cap C = \{a, b\}$. Then $\text{diam}(C) \leq 5\lambda(E)$.*

**Proof.** Let $D$ be the interior of the simple closed curve $C \cup E$. Obviously $D \subset G$. Let $p \in C$ and $q \in E$ satisfy
$$|p - q| = \sup\{|z - w| : z \in C \text{ and } w \in E\}.$$
Let $C'$ be the elementary $j$-characteristic containing $p$. Since a $j$-characteristic can have at most one point in common with an $i$-characteristic, there is a point $q' \in E$ such that $C'$ contains a subarc $J = pq'$ whose interior lies in $D$. Let $\text{diam}(C) = |z_1 - z_2|$, where $z_1$, $z_2 \in C$. Then
$$\text{diam}(C) = |z_1 - z_2| \leq |z_1 - a| + |a - b| + |b - z_2| \leq 2|p - q| + \lambda(E),$$
since $|z_1 - a|, |b - z_2| \leq |p - q|$. But
$$|p - q| \leq |p - q'| + |q' - q| \leq \lambda(J) + \lambda(E),$$
so that
$$\text{diam}(C) \leq 2\lambda(J) + 3\lambda(E).$$

But it follows by a simple argument based on Proposition 2.7 (Length Monotonicity) that $\lambda(J) \leq \lambda(E)$, since the elementary $i$-characteristic through each interior point $z$ of $J$ contains a subarc joining two distinct points of $E$. Thus $\text{diam}(C) \leq 5\lambda(E)$ as claimed. ∎

**Proposition 3.7.** *Let $C_1 = pq_1$ and $C_2 = pq_2$ be distinct elementary $i$-characteristics of a $K$-quasi-HP function on $G$ with an endpoint $p \in \partial G$ in common. Let $J \subset G$ be a $j$-characteristic arc joining a point of $c_1$ of $C_1 \cap G$ to a point $c_2$ of $C_2 \cap G$, so that $C_1$, $C_2$ and $J$ form the three sides of a "characteristic triangle" $T$. Let $P \subset J$ denote the set of points at which $J$ is not strictly concave towards the inside of $T$. Then $\lambda(P) = 0$.*

**Proof.** Suppose not. Then since the elementary $i$-characteristic passing through any point of $J$ must exit at $p$, after replacing the original $J$ by an appropriate subarc and changing $C_1$ and $C_2$ accordingly we can assume that $\lambda(P \cap J) = \epsilon > 0$, $\lambda(\theta(J)) < \frac{1}{100K}$ and $\lambda(N) < \frac{\epsilon}{8}$, where $N = J \setminus P$. Let $C_k$ be parametrized by $z_k(s)$, $0 \leq s \leq \lambda_k$, with $z(0) = p$. Let $\{k, l\} = \{1, 2\}$. Let $J_k(s)$ denote the $j$-arc joining $z_k(s)$ to a point $w_k(s) \in C_l \cap G$. Note that we only know that $J_k(s)$ is defined for $s \leq \lambda_k$ sufficiently near $\lambda_k$. It follows from length monotonicity (Proposition 2.7) that $\lambda(J_k(s)) \geq \lambda(P_k(s)) \geq \epsilon$, where $P_k(s)$ is the set of points of $J_k(s)$ joined to points of $P$ by an $i$-arc. By the quasi-HP property $\lambda(\theta(J_k(s))) < \frac{1}{100}$, so that $J_k(s)$ is almost straight and in particular the distance between its endpoints is at least $\frac{1}{2}\lambda(J_k(s)) \geq \frac{1}{2}\lambda(P_k(s)) \geq \frac{\epsilon}{2}$. Let $\xi_k$ be the infimum of all $s$ for which $J_k(s)$ is defined. Since the distance between the endpoints of $J_k(s)$ is at least $\frac{\epsilon}{2}$, it is clear that at least one of $\xi_1$, $\xi_2$ must be positive, and for definiteness we assume that $\xi_1 > 0$. For $\sigma \in (\xi_1, \lambda_1]$ let $J_1(\sigma)$ be parametrized by $\zeta(s, \sigma)$, $0 \leq s \leq \lambda(J_1(\sigma))$, with $\zeta(0, \sigma) \in C_1$. It is clear that there are $\delta, T > 0$ such that

(3.2) $\qquad \text{dist}(\zeta(s, \sigma), C_1 \cup C_2) \geq Ts$, for $s \in (0, \delta)$, $\sigma \in (\xi_1, \lambda_1]$.

From the fact that a $j$-arc can intersect an $i$-arc at most once in $G$ it easily follows that for each point $z \in C_1 \cap G$ there is a $\delta_1 = \delta_1(z)$ such that

(3.3) $\qquad |z - \zeta(s, \sigma)| \geq \delta_1$, for $s \in [\delta, \lambda(J_1(\sigma))]$, $\sigma \in (\xi_1, \lambda_1]$.

From (3.2) and (3.3) together with the fact that for $\sigma \in (\xi_1, \lambda_1]$, $J_1(\sigma) = J_2(\sigma')$ for some $\sigma' \in (\xi_2, \lambda_2]$ it follows that as $\sigma \to \xi_1$, $J_1(\sigma)$ tends to an arc $J_0$ which contains $p$, which joins $z_1(\xi_1)$ to the point $z_2(\xi_2)$ of $C_2$ and the distance between whose endpoints is at least $\frac{\epsilon}{2}$. Furthermore, either $J_0$ consists of a $j$-arc joining $z_1(\xi_1)$ to $p$ or (in the case that $\xi_2 > 0$) it consists of such an arc together with another $j$-arc joining $p$ to $z_2(\xi_2)$. One of the endpoints, which we henceforth call $q$, is at a distance of at least $\frac{\epsilon}{4}$ from $p$. For definiteness or by renaming we assume that $q = z_1(\xi_1) \in C_1$. Let $J'_0$ be the subarc of $J_0$ which joins $q$ to $p$. Let $A$ denote the arc $pq$ of $C_1$. Then $J'_0$ and $A$ form the two sides of a "characteristic bilateral" $B$. Let $P'_0$ be the subset of points of $J'_0$ which correspond to (i.e., are joined by $i$-arcs to) points of $P$, (that is, $P'_0$ is the set of points at which $J_0$ is nonconcave towards the inside of $B$), and let $N'_0 \subset J'_0$ be the subset corresponding to $N$. Since by length monotonicity, $\lambda(N'_0) \leq \lambda(N) \leq \frac{\epsilon}{8}$, it follows that

(3.4) $\qquad \lambda(P'_0) \geq \lambda(J'_0) - \lambda(N'_0) \geq \frac{\epsilon}{4} - \frac{\epsilon}{8} = \frac{\epsilon}{8}$.

Let $J'_0$ be parametrized by $z_0(s)$, $0 \leq s \leq \lambda(J'_0)$, with $z_0(0) = p$, and for $s \in (0, \lambda(J'_0))$ let $A(s)$ be the part of the elementary $i$-characteristic through $z_0(s)$ in $\overline{B}$, so that $A(s)$ joins $z_0(s)$ to $p$ (since its interior can cross neither $A$ nor $J'_0$). It follows from Proposition 3.6 that $\text{diam}(A(s))$ tends to 0 as $s \to 0$. Let $qp'$ be an initial arc of $A$ for which $\lambda(\theta(qp')) < \frac{1}{100K}$. Since $\lambda(\theta(W)) < \frac{1}{100}$ for any translate $W$ of either $J_0'$ or $qp'$ by the quasi-HP-property, it follows that we are able to translate $qp'$ in $B$ all the way down $J_0' \setminus \{p\}$ from $q$ without meeting the boundary of $B$, since for simple geometric reasons all these translates, being essentially perpendicular to the virtually straight arc $J_0'$ must stay away from $p$. If $A'(s)$ denotes the translate of $qp'$ with initial point $z_0(s)$, then $A'(s) \subset A(s)$, and therefore $\text{diam}(A'(s)) \to 0$ as $s \to \lambda(J_0')$. This means that each point $z \in qp'$ is joined to $p$ in $B$ by a $j$-arc $J''(z)$ such that if $P''(z)$ is the set of points of $J''(z)$ joined to points of $P$ by $i$ arcs in $T$, then by (3.4) and length monotonicity $\lambda(P''(z)) \geq \frac{\epsilon}{8}$ and by the quasi-HP-property $\lambda(\theta(J''(z))) < \frac{1}{100}$, so that
$$(3.5) \qquad |z - p| \geq \frac{\epsilon}{16}$$
We can now repeat this process starting with $J''(p')$ instead of $J_0' = J''(q)$, and continue doing so to obtain in the end $J''(z)$ for all $z \in A \setminus \{p\}$. However, by the argument we just gave we now have (3.5) for all $z \in A$, which is absurd since $p$ is an endpoint of $A$. Therefore $\lambda(P) = 0$.∎

**Proposition 3.8.** *The two endpoints of a nontrivial elementary characteristic of a quasi-HP function on $G$ must be different.*

**Proof**. This follows easily from Proposition 3.6∎

**Proposition 3.9.** *Two elementary $i$-characteristics of a quasi-HP function on $G$ with the same endpoints must be identical.*

**Proof**. Assume to the contrary that points $p_1 \neq p_2$ of $\partial G$ are joined by distinct elementary $i$-characteristics $C_1$ and $C_2$. It follows from Proposition 3.8 that the elementary $i$-characteristic through any point of the simply connected domain $D$ bounded by $C_1 \cup C_2$ must also have endpoints $p_1$ and $p_2$. But then it follows easily from Proposition 3.7 that all $j$-characteristic arcs in $D$ are straight line segments. But this contradicts Proposition 3.7.∎

**Proposition 3.10.** *Two different elementary characteristics of a quasi-HP function on $G$ which both exit $G$ regularly at $p \in G$ cannot be tangent to each other there.*

**Proof**. This is an easy consequence of Proposition 3.7 and the quasi-HP property.∎

**Definition 3.11.** Let $\theta$ be a quasi-HP function on $G$ and let $C_0$ be an elementary $i$-characteristic of $\theta$ with endpoints $a, b \in \partial G$. Let $B$ be one of the boundary arcs of $\partial G$ with endpoints $a, b$. Then the subdomain $D$ of $G$ for which $\partial D = C_0 \cup B$ will be called an *$i$-characteristic subdomain*.

When we wish to indicate the elementary $i$-characteristic involved, we will denote the $i$-characteristic subdomain by $(D, C_0)$. The arc $B = \overline{\partial D - C_0} \subset \partial G$, called the *bottom* of $(D, C_0)$ and denoted by $\text{bot}(D)$, will be considered to have the order "$<$" corresponding to the positive orientation of $\partial D$. We shall freely use interval notation as well as the terms, "to the right of", "to the left of", "between", etc. when dealing with $\text{bot}(D)$. Furthermore if $ab$ is an elementary $i$-characteristic joining points $a$, $b$ of $\text{bot}(D)$, it will be understood that $a \leq b$, unless otherwise indicated. When dealing with a characteristic subdomain $(D, C_0)$, we shall work with the class $\mathcal{I}(D)$ of nontrivial elementary $i$-characteristics $C$ which join points $p, q$ of $\text{bot}(D)$. Note that $C_0 \in \mathcal{I}(D)$. If $C = pq \in \mathcal{I}(D)$, then "above" $C$ refers to the part of $D$ not in the closed region bounded by $C \cup [p, q]$. If $F \subset \overline{D}$ is a compact set for which $F \cap \text{bot}(D) \neq \emptyset$, we say that an elementary $i$-characteristic $E = pq \in \mathcal{I}(D)$ *envelopes* $F$, and write $F \preceq E$, if $F$ is contained the closed set bounded by the simple closed curve $[p, q] \cup pq$. We only deal with sets $F$ for which each component of $F$ has points in $\operatorname{bot}(D)$. For two such sets $F_1$, $F_2 \subset \overline{D}$ we say that $F_1$ *precedes* $F_2$ and write $F_1 \leq F_2$ if for all, $f_1$, $f_2$ with $f_k \in F_k \cap \operatorname{bot}(D)$, $k = 1, 2$ there holds $f_1 \leq f_2$. It is clear that if $C_1, C_2 \in \mathcal{I}(D)$, then one of the following is true: $C_1 \preceq C_2$, $C_1 \leq C_2$ or $C_2 \leq C_1$. The first of these possibilities includes the case $C_1 = C_2$.

In what follows $L(G, p, d)$ will denote the segment $\overline{pq}$ of length $d$ which is orthogonal to $\partial G$ at $p$, which emanates from $p$ into $G$ and which is oriented from $p$ to $q$.

**Proposition 3.12.** *There are an absolute constant $\omega_0 = \omega_0 < 2$ and a constant $\omega_1 = \omega_1(K, \rho)$ with the following properties. Let $G \in \mathcal{G}(\rho)$ and $\theta \in \operatorname{HP}(G, K)$, and let $C$, parametrized by $z(s)$, $0 \leq s \leq \lambda(C)$, be an arc of an i-characteristic the distance between whose endpoints is at least $\rho$, for which*
(3.6) $$C \cap L(G, p, \omega_0 \rho) = \{z(0)\},$$
*which emanates to the right (left) of $L(G, p, \omega_0 \rho)$ and for which each j-half-characteristic emanating to the left (right) of $C$ has length at least $\rho$. Let $\partial G$ be parametrized by $w(s)$, $0 \leq s < \lambda(\partial G)$, respectively with $z(0)$, $w(0) \in L(G, p, \omega_0 \rho)$ and $|\arg\{z'(0)\} - \arg\{w'(0)\}| \leq \frac{\pi}{2}$.*

*Then*
$$\arg\{z'(0)\} \geq \arg\{w'(0)\} - \omega_1 \sqrt{|z(0) - p|},$$
*if $C$ emanates to the right of $L(G, p, \xi_0)$, and*
$$\arg\{z'(0)\} \geq \arg\{w'(0)\} + \omega_1 \sqrt{|z(0) - p|},$$
*if $C$ emanates to the left of $L(G, p, \omega_0)$.*

Proof. Clearly it is enough to handle the case in which $C$ emanates to the right of $L(G, p, \rho)$. Without loss of generality we can assume that $p = 0 = \arg\{w'(0)\}$, so that in particular $D = N(-\rho i, 0) \subset \mathbb{C} \setminus G$ and $\partial N(-\rho i, 0)$ is tangent to $\partial G$ at 0. Let $R_d = \{di + te^{-i\alpha} : t \geq 0\}$, $\alpha = \alpha(d) \in (0, \frac{\pi}{2})$ be the ray emanating to the right of $L(G, 0, \rho)$ from the point $di \in L(G, 0, \rho)$ which is tangent to $\partial D$; and let the point of tangency be $z_d$. Then
$$\cos\alpha = \frac{\rho}{\rho + d} = 1 - \frac{d}{\rho} + O(\frac{d^2}{\rho^2}),$$
so that there is some $A_1 < 2$ such
(3.7) $$2\sqrt{\frac{d}{\rho}} \geq \alpha \geq \frac{1}{2}\sqrt{\frac{d}{\rho}}, \text{ for } d \leq A_1 \rho.$$
Let $T_d$ be the curvilinear triangle bounded by $L(G, 0, d)$, the line segment $[di, z_d]$ and the (shorter) arc of $\partial D$ with endpoints 0, $z_d$. Then it is easy to see that
(3.8) $$\operatorname{diam}(T_d) \leq \rho\alpha + d \leq 2\rho\sqrt{\frac{d}{\rho}} + d = 2\sqrt{\rho d} + d, \text{ for } d \leq A_1\rho,$$
so that (after replacing $A_1$ by a smaller value, if necessary)
(3.9) $$\operatorname{diam}(T_d) < \rho, \text{ for } d \leq A_1\rho.$$
Now assume that $C$ is as in the hypothesis with $z(0) = di$, where $d \leq A_1\rho$ and that
$$-\beta = \arg\{z'(0)\} \leq -N\alpha.$$
Here $N$ is some (large) number yet to be determined. Then, since by (3.9) $\operatorname{diam}(C) \geq \rho > \operatorname{diam}(T_d)$ and because obvious angle considerations imply that $C$ enters $T_d$ at $di$, in light of assumption (3.6), $C$ must exit $T_d$ at a point of the segment $[di, z_d]$. Say $s_0$ is the smallest value of $s$ for which $z(s) \in [di, z_d]$. Then obviously $\theta([0, s_0]) \supset [-\beta, -\alpha]$, so that $[0, s_0]$ has a subinterval $[s_1, s_2]$ for which $\theta([s_1, s_2]) = [-\beta, -\alpha]$. Since $\beta - \alpha \leq \frac{\pi}{2}$ it follows that
$$s_2 - s_1 = \lambda(z([s_1, s_2])) \leq 2|z(s_2) - z(s_1)| \leq 2\operatorname{diam}(T_d) \leq 2(2\sqrt{\rho d} + d),$$
by (3.8). Thus by the mean value theorem there must be some $s \in [s_1, s_2]$ for which
$$(3.10) \qquad D^+\theta(z(s_0)) \geq \frac{(N-1)\alpha}{2(2\sqrt{\rho d}+d)}.$$
But by the curvature bound (Proposition 2.6) and the hypothesis regarding the $j$-characteristics $D^+\theta(z(s_0)) \leq \frac{K}{\rho}$, so that in light of (3.10) and the lower bound in (3.7) we have
$$\rho(N-1)\sqrt{\tfrac{d}{\rho}}/(8\sqrt{\rho d} + 4d) \leq K.$$
But the left hand side of this inequality is $(N-1)/(8 + 4\sqrt{\tfrac{d}{\rho}}) \geq \frac{N-1}{8+4\sqrt{A_1}}$ for $d \leq A_1\rho$, so that we have a contradiction for $N = (8 + 4\sqrt{A_1})K + 2$. Thus if $d \leq A_1\rho$ we must have
$$\arg\{z'(0)\} \geq -2\frac{(8+4\sqrt{A_1})K+2}{\sqrt{\rho}}\sqrt{d},$$
in light the upper bound in (3.7). Thus we have proved the proposition with $\omega_0 = A_1$ and $\omega_1 = 2\frac{(8+2\sqrt{A_1})K+2}{\sqrt{\rho}}$. ∎

**Proposition 3.13.** *There are positive constants $\xi_0 = \xi_0(K, \rho)$ and $\xi_1 = \xi_1(K, \rho)$ with the following properties. Let $G \in \mathcal{G}(\rho)$, $\theta \in \operatorname{HP}(G, K)$, and let $C = ab$ be an elementary $i$-characteristic, with corresponding characteristic subdomain $(D, C)$. Let $C$ and $A = \operatorname{bot}(D)$ be given by $z(s)$, $0 \leq s \leq \lambda(C)$ and $w(s)$, $0 \leq s \leq \lambda(A)$, respectively, with $z(0) = w(0) = a$, $a$ being the leftmost point of $\operatorname{bot}(D)$. Then for each $p = w(\sigma_0) \in A$ for which $\operatorname{dist}(p, \partial G \setminus A) \geq \rho$, the segment $L(G, p, \xi_0)$ contains at most one point of $C$, and if there is such a point $z(s)$, then there holds*
$$(3.11) \qquad |\arg\{z'(s)\} - \arg\{w'(\sigma_0)\}| \leq \xi_1\sqrt{|z(s) - p|}.$$

**Proof.** Without loss of generality we may assume that $p = 0 = \arg\{w'(\sigma_0)\}$. Let $\omega_0$ and $\omega_1$ be as in the preceding proposition. Let $L(\epsilon) = L(G, p, \epsilon)$. We assume for the moment that $z(s) \in L(\omega_0 \frac{\rho}{2})$ is the only point of $C$ on this segment. Since immediately to the left of $C$ there are points in the complement of $D$ and the interior of the segment $L(|z(s)|)$ lies in $D$, it is clear that if we write $\arg\{z'(s_0)\} = e^{i\tau}$, with $-\pi < \tau \leq \pi$, then $|\tau| \leq \frac{\pi}{2}$. Let $E$ be one of the two subarcs of $C$ one of whose endpoints is $z(s)$ and the other of which is a point $q \in G$ for which $|z(s) - q| = \frac{\rho}{2}$. It follows from the preceding proposition that if $\omega_1(K, \frac{\rho}{2})\sqrt{|z(s)|} \leq \frac{\pi}{4}$, then $C$ cannot be tangent to $L(\omega_0 \frac{\rho}{2})$. Thus, of the two arcs $E$ one moves to the right of $L(\omega_0 \frac{\rho}{2})$ as we move along it away from $z(s)$ and the other moves to the left. But then by the preceding proposition we have (3.11) with
$$(3.12) \qquad \xi_1 = \omega_1(K, \tfrac{\rho}{2})$$
for any $z(s)$ for which
$$(3.13) \quad |z(s)| \leq \min\{\omega_0 \tfrac{\rho}{2}, \left(\tfrac{\pi}{4\omega_1(K,\rho/2)}\right)^2\} \leq \min\{\omega_0 \tfrac{\rho}{2}, \left(\tfrac{\pi}{4\omega_1(K,\rho/2)}\right)^2, \tfrac{\pi\rho}{21K}\} = \xi_0.$$
Assume $\epsilon \leq \xi_0$ and that $L(\epsilon)$ contains at least two points of $C$. Then there are $s'$, $s'' \in (0, \lambda(C))$, $s' < s''$, such that the interior of $L(\min\{|z(s')|, |z(s'')|\})$ contains no point of $C$ (and is therefore contained in $D$) and $z(s')z(s'') \cap L(\epsilon) = \{z(s'), z(s'')\}$. By Proposition 3.6, $\operatorname{diam}(z(s')z(s'')) \leq 5\epsilon < \omega_0 \frac{\rho}{2}$. There are the following two cases.

(i) $|z(s')| < |z(s'')|$. In this case, simple topological arguments show that $z(s')z(s'')$ must lie to the right of $L(\epsilon)$. Also, by the foregoing and our definition of $\xi_0$, one easily has
$$|\arg\{z'(s')\} - \arg\{w'(\sigma_0)\}| \leq \tfrac{\pi}{4}.$$
But then since $C$ crosses $L(\epsilon)$ again at $z(s'')$, we must have that for some $t' < t''$ in $(s', s'')$, $\arg\{z'(t'')\} - \arg\{z'(t')\} \geq \tfrac{\pi}{2} - \tfrac{\pi}{4} = \tfrac{\pi}{4}$, so that by mean value considerations as in the proof of the preceding proposition together with the curvature bound we see that there must be a point $\sigma' \in (t', t'')$ at such that
$$\tfrac{K}{\rho/2} \geq \tfrac{d\arg\{z'(s)\}}{ds}|_{s=\sigma'} \geq \tfrac{\pi/4}{2\mathrm{diam}(z(s')z(s''))} \geq \tfrac{\pi}{10\epsilon} \geq \tfrac{\pi}{10\xi_0},$$
which implies that $\xi_0 \geq \tfrac{\pi\rho}{20K}$, a contradiction, since $\xi_0 \leq \tfrac{\pi\rho}{21K}$.

(ii) $|z(s')| > |z(s'')|$. In this case, simple topological arguments show that $z(s')z(s'')$ must lie to the left of $L(\epsilon)$ and we proceed analogously to the way we did in case (i). This completes the proof of the proposition with $\xi_1$ and $\xi_0$ as defined in (3.12) and (3.13), respectively. ■

We now extend of the notion of characteristic to include certain arcs whose interiors contain points of $\partial G$. Although what follows does not give the most exhaustive extension possible, it is sufficient for our present needs. Let $(D_0, C_0)$ be an $i$-characteristic subdomain of $G$ and consider a monotone decreasing sequence $\{C_k = a_k b_k : k \geq 0\}$ in $\mathcal{I}(D_0)$. In other words, $a_k \leq a_{k+1} < b_{k+1} \leq b_k$, $k \geq 0$, so that the arcs $A_k = [a_k, b_k]$ of $\mathrm{bot}(D_0) = A_0$ are nested. Let $D_k \subset D_0$ be the $i$-characteristic subdomain bounded by $C_k \cup A_k$, so that $D_k \supset D_{k+1}$, $k \geq 0$. We regard $C_k$ as being oriented from $a_k$ to $b_k$. Let $a_k \to a$ and $b_k \to b$. We define $C$ to be the set of limits of sequences $\{z_l\}$, where $z_l \in C_{k_l}$, $k_l \to \infty$. Any such $C$ will be called an *extended characteristic with endpoints $a$ and $b$*. An extended characteristic consisting of a single point $p \in \partial G$ will be called *trivial*. The following proposition, in which the notation is the same as that of the immediately preceding sentences, contains the basic properties of extended characteristics. Note that $M(G, K)$ is, as in Proposition 3.4, an upper bound on the length of elementary characteristics of $K$-quasi-HP functions on $G$.

**Proposition 3.14.** *Let $C$ be an extended characteristic with endpoints $a$ and $b$. If $a = b$ then $C = \{a\}$. Otherwise, $C$ is a simple arc joining $a$ to $b$ and $\lambda(C) \leq M(G, K)$. If $C$ is parametrized by $z(s)$, $0 \leq s \leq \lambda(C)$ with $z(0) = a$, then $C \cap \partial G \subset [a, b]$ and for points of $C \cap [a, b]$ the order with respect to $\mathrm{bot}(D_0)$ coincides with the order with respect to $s$. Furthermore, the function $z$ is continuously differentiable on $(0, \lambda(C))$ and for $t \in (0, \lambda(C))$, $z(t)$ is joined to a point of $\partial G \backslash (a, b)$ by a unique $j$-characteristic arc $J(t)$ emanating to the left of $C$ and*

(3.14) $$D^+\arg\{z'(t)\} \leq \tfrac{K}{\lambda(J(t))}.$$

**Proof**. We begin by observing that

(3.15) $$A_0 \cap C \subset [a, b].$$

To see this, note that for each $w \in A_0 \backslash [a, b]$ there is an $n$ for which $w \notin A_n$. Since $C_n \cap A_0 = \{a_n, b_n\}$, $w$ is not in $C_n$ either, and therefore, $w \notin \overline{D_n}$. But then, by the monotonicity of $\{\overline{D_k}\}$,
$$\mathrm{dist}(w, C_k) \geq \mathrm{dist}(w, \overline{D_k}) \geq \mathrm{dist}(w, \overline{D_n}) > 0 \text{ for } k \geq n,$$
from which the desired conclusion follows at once. From (3.15) it follows immediately that $C \cap \partial G \subset [a, b]$.

Next we note that if $a = b$, then Proposition 3.6 implies that $\mathrm{diam}(C_k) \to 0$, so that $C = \{a\}$. For the remainder of the proof we therefore assume that $a < b$.

By this assumption and Proposition 3.4 on the boundedness of the lengths of characteristics there exist $l_1$ and $l_2$ such that $0 < l_1 \leq \lambda(C_k) \leq l_2 < \infty$. Let $C_k$ be parametrized by $z = z_k(s)$, $0 \leq s \leq \lambda(C_k)$. Let $n_0$ be such that
$$|a_k - a|, |b_k - b| < |b - a|/3 \text{ for } k \geq n_0.$$
Let $0 < \epsilon < |b - a|/3$, and for $k \geq n_0$ let
$$\alpha_k = \sup\{s : z_k(s) \in \partial N(a, \epsilon)\} \quad \text{and} \quad \beta_k = \inf\{s : z_k(s) \in \partial N(b, \epsilon)\},$$
and let $E_k = E_k(\epsilon) = z_k([\alpha_k, \beta_k])$. Then
$$(3.16) \qquad \delta = \inf\{\mathrm{dist}(E_k, \partial G \backslash A_k) : k \geq n_0\} > 0,$$
since otherwise there would be a point of $C$ in $A_0 \backslash [a, b]$, in contradiction of (3.15). For $k \geq n_0$ and $s \in [\alpha_k, \beta_k]$, let $J_k(s)$ denote the $j$-half-characteristic emanating to the left of $C_k$ and joining $z_k(s)$ to a point $w_k(s)$ of $\partial G \backslash A_k$. In light of (3.16),
$$(3.17) \qquad \lambda(J_k(s)) \geq \mathrm{dist}(z_k(s), \partial G \backslash A_k) \geq \delta > 0 \text{ for all } k \geq n_0, s \in [\alpha_k, \beta_k].$$
By the curvature bound this means that for each $\epsilon > 0$ we have an upper bound on the curvature to the left of $E_k$. More precisely, we have that
$$(3.18) \qquad D^+\theta(z_k(s)) \leq \frac{K}{\lambda(J_k(s))} \leq \frac{K}{\mathrm{dist}(z_k(s), \partial G \backslash A_k)} \leq \frac{K}{\delta}, \text{ for } k \geq n_0, s \in [\alpha_k, \beta_k].$$
The curvature bound trivially implies that
$$(3.19) \qquad D^-\theta(z_k(s)) \geq -\frac{K}{\mathrm{dist}(z_k(s), \partial G)}.$$
There is a neighborhood $U$ of $\partial G$ such that for each $z \in U$ there is a unique $p(z) \in \partial G$ for which $|z - p(z)| = \min\{\mathrm{dist}(z, \zeta) : \zeta \in \partial G\}$ and for which $p$ is continuous. For $q$ in $\partial G$, let $e^{i\phi(q)}$ be the positively oriented unit tangent to $\partial G$ at $q$, so that in $U$, $e^{i\phi(p(z))}$ is continuous. It follows from (3.18) and Proposition 3.10 that there is an $\Lambda = \Lambda(\epsilon)$ such that
$$(3.20) \qquad |z_k'(s) - e^{i\phi(p(z_k(s)))}| \leq \Lambda \sqrt{|z_k(s) - p(z_k(s))|}.$$
Bounds (3.18), (3.19), (3.20) and the imply that the family $\{z_k'(s)\}$ is uniformly bounded and equicontinuous (that is, has a uniform modulus of continuity valid for all $k \geq n_0$ on the $[\alpha_k, \beta_k]$). From this it follows that there are $\alpha$, $\beta$ and a sequence $\{k_l\}$ such that $\alpha_{k_l} \to \alpha$ and $\beta_{k_l} \to \beta$ and $z_{k_l}(s) \to z_\epsilon(s)$ uniformly (in the obvious sense), where $z_\epsilon$ is continuously differentiable and parametrizes an arc $C(\epsilon)$ which joins point of $\partial N(a, \epsilon)$ to point of $\partial N(b, \epsilon)$ in the part of $\overline{D_0}$ lying outside of both these circles. The arc $C(\epsilon)$ is simple since by Proposition 3.13 any sufficiently small neighborhood of a point of $(a, b)$ can contain at most a single arc of $C_k$. Clearly, $\lambda(C(\epsilon)) \leq M(G, K)$. The nested nature of the $D_k$ then implies that (a least for sufficiently small $\epsilon$) the entire sequence converges to $C(\epsilon)$. Also, it is clear that $C(\epsilon')$ is an extension of $C(\epsilon)$, for $\epsilon' < \epsilon$. We have that $C = \cup \{C(\epsilon) : \epsilon > 0\}$, since $\mathrm{diam}(z_k([0, \alpha_k]))$ and $\mathrm{diam}(z_k([\beta_k, \lambda(C_k)]))$ tend to 0 uniformly in $k$ as $\epsilon \to 0$ by Proposition 3.6. It follows from $C = \cup \{C(\epsilon) : \epsilon > 0\}$ that $\lambda(C(\epsilon)) \leq M(G, K)$ and that $C$ is a simple arc with endpoints $a$, $b$ and furthermore that $C$ is parametrized by a function $z(s)$ which is continuously differentiable on $(0, \lambda(C))$. That the two possible orderings of the points of $C \cap [a, b]$ coincide follows from the fact that $C$ is a simple arc and (3.15).

The existence and uniqueness of $J(t)$ is trivial when $z(t) \in G$. When $z(t) \in (a, b)$ the existence of $J(t)$ follows from a straightforward compactness argument. In light of (3.20), for such $t$, any corresponding $J(t)$ must be orthogonal to $\partial G$ at $z(t)$, so that the uniqueness of $J(t)$ follows from Proposition 3.10. Bound (3.14) follows from (3.18). ∎

We will refer to an extended characteristic $C$ joining $a, b \in \partial G$ as $ab$; points of $ab \cap \partial G$ will be called *contact points* and contact points other than $a$ and $b$ will be called *proper contact points*. It is clear that $\theta$ can be continuously extended to $G \cup (ab \setminus \{a, b\})$. For what is to follow it is important to understand that if $(D_0, C_0)$ is an $i$-characteristic subdomain, then in addition to the extended characteristics constructed above, $D_0$ might contain extended characteristics $C'$ joining points $a' < b'$ arising from a sequence of $i$-characteristic subdomains $\{(D'_k, C'_k)\}$, where $C'_k = a'_k b'_k$ with $a'_{k+1} \leq a'_k < b'_k \leq b'_{k+1}$ (where the order is with respect to $A_0 = \text{bot}(D_0)$, as above). Here the $C'_k$ are contained in $D_0$ but the other part of the boundary of $D'_k$ is $\partial G \setminus (a'_k, b'_k)$ (where here again the interval notation refers to the order on $A_0 = \text{bot}(D_0)$). Note that if for such an extended characteristic $C'$, $C' \cap \partial G$ has points other than $a'$ and $b'$ (that is, if $C'$ is not simply an elementary characteristic), then the contact points will not occur monotonically with respect to the order on $A_0$ when $C'$ is traversed from $a'$ to $b'$. *With respect to the characteristic subdomain $D_0$ the extended characteristics $C$ constructed originally will be referred to as monotone*, and this other kind of extrended characteristic $C'$ with proper contact points will be said to be *nonmonotone*. We consider that an extended characteristic $ab$ exits $G$ at all of its contact points, and we use the terms "exits regularly" and "exits singularly" at $a$ or $b$ in the obvious fashion. Clearly, $ab$ exits regularly at its proper contact points.

**Proposition 3.15.** *Let $C$ be an extended $j$-characteristic which exits at $p$ and is parametrized by $z = z(s)$, $0 \leq s \leq L$, with $z(0) = p$. Let $\phi(s) = \arg\{z'(s)\}$ and assume that $\lim_{s \to 0} \phi(s)$ does not exist. For $s \in (0, L)$ let $E(s)$ be the elementary $i$-characteristic containing $z(s)$ and for $\epsilon, T > 0$ let*
$$P(\epsilon, T) = \{s < \epsilon : \phi'(s) \geq T\} \quad \text{and} \quad N(\epsilon, T) = \{s < \epsilon : \phi'(s) \leq -T\}.$$
*Then for all $\epsilon, T > 0$, $P(\epsilon, T)$ and $N(\epsilon, T)$ have positive measure and $\text{diam}(E(s)) \to 0$ as $s \to 0$.*

**Proof**. To see that $\lambda(P(\epsilon, T)) > 0$, assume to the contrary for some $\epsilon_0, T_0 > 0$, $\lambda(P(\epsilon_0, T_0)) = 0$, so that $\phi'(s) \leq T_0$ a.e. on $(0, \epsilon_0]$. Then $\phi = \phi^+ + \phi^-$ on $(0, \epsilon_0]$, where $\phi^+$ is Lipschitz continuous and nondecreasing and $\phi^-$ is continuous and nonincreasing. From this in turn it follows that either $\phi^-$ has a finite limit as $s \to 0$ or it tends to $\infty$ as $s \to 0$, so that in fact the latter is the case since $C$ exits singularly. But this means that $C$ spirals around $p$, which is clearly impossible. Thus $P(\epsilon, T)$ must have positive measure for all $\epsilon, T > 0$. One sees similarly that $N(\epsilon, T)$ has positive measure. That $\text{diam}(E(s)) \to 0$ follows immediately from the characteristic length bound. ∎

Let $(D, C_0)$ be a characteristic subdomain of $G$. If $F \subset D \cup \text{bot}(D)$ is a compact set for which $F \cap \text{bot}(D) \neq \emptyset$, then the family $\mathcal{C}(F) = \{C \in \mathcal{I}(D) : F \preceq C\}$ is clearly linearly ordered with respect to the relation $\preceq$. From this it easily follows that there is a unique monotone extended $i$-characteristic $E_0 = ab$ for which $F \cap \text{bot}(D) \subset [a, b]$ and $E_0 \preceq C$ for all $C \in \mathcal{C}(F)$.

**Notation 3.16.** This monotone extended $i$-characteristic $E_0$ will be denoted by $\min_D(F)$. In the statement and proof of the following proposition all order relations are with respect to $(D, C_0)$.

**Proposition 3.17.** (Structure of $\min_D(A \cup B)$) *Let $(D, C_0)$ be an $i$-characteristic subdomain of $G$. Let $a_1 < a_2 < b_1 < b_2$ be such that $a_1 a_2 = A$ and $b_1 b_2 = B$ belong to $\mathcal{I}(D)$ and let $U = \min_D(A \cup B) = ef$, $e < f$. Then least one of the following things must happen.*
(*i*) *$U$ has a contact point between $A$ and $B$,*
*or $ef \in \mathcal{C}(A \cup B)$ and there is a nonmonotone extended $i$-characteristic $ab$ with $e \le a \le b \le f$ and at least one of the following happens:*
(*ii*) *$ea$ contains a subarc $L \in \mathcal{I}(D)$ for which $A \preceq L \preceq U$, but $B \preceq L$ doesn't hold.*
(*iii*) *$bf$ contains a subarc $L \in \mathcal{I}(D)$ for which $B \preceq L \preceq U$, but $A \preceq L$ doesn't hold.*

**Proof.** It follows from the hypotheses that $e < f$ (i.e., that $e \ne f$). We show that if (*i*) does not hold, then at least one of (*ii*) or (*iii*) does. Let
$$e' = \sup\{z : z \in \text{bot}(D) \cap \min_D(A \cup B), z \le a_1\}$$
and
$$f' = \inf\{z : z \in \text{bot}(D) \cap \min_D(A \cup B), z \ge b_2\}.$$
Since (*i*) doesn't hold it is easy to see that the subarc $e'f'$ of $ef$ is a member of $\mathcal{C}(A \cup B)$, so that in fact $e'f' = ef$. Obviously, $ef \cap (A \cup B) \subset \{e, f\}$, since otherwise $ef$ would have to be $A$ or $B$. Let $w_0$ be any point of $ef$ other than $e$ or $f$, and consider a small $j$-arc $E$ whose initial point is $w_0$, which extends to the right of $ef$ (that is, into the simply connected domain bounded by $ef \cup [e, f]$) and which is disjoint from $A \cup B$. Let $E$ be parametrized by $z = w(s)$, $0 \le s \le \lambda_0$ with $w(\lambda_0) = w_0$, and let $C(s)$ denote the elementary $i$-characteristic through $w(s)$. Note that
(3.21) $$C(s_1) \preceq C(s_2) \text{ for } s_1 \le s_2.$$
Let the left and right endpoints of $C(s)$ be $l(s)$ and $r(s)$. By the minimality of $ef$ it must be that for no $s \in (0, \lambda_0)$ can we have both $l(s) \le a_1$ and $r(s) \ge b_2$. From (3.21) it therefore follows that either for all $s \in (0, \lambda_0)$, $l(s) > a_1$, or for all $s \in (0, \lambda_0)$, $r(s) < b_2$. Assume the latter occurs. Then in fact $r(s) \le b_1$ for all $s \in (0, \lambda_0)$ since otherwise $C(s)$ would cross $B$ in $D$, which is impossible since they are distinct elementary $i$-characteristics. Consider the $i$-characteristic subdomain $(D', C(0))$ bounded by $C(0)$ and $\partial G \setminus (l(0), r(0))$. If we take any sequence $\{s_k\}$ in $(0, \lambda_0)$ which tends monotonically to $\lambda_0$, and consider the corresponding sequence of elementary $i$-characteristics $C_k = C(s_k)$ in $D'$, it is clear that they will give rise to an extended $i$-characteristic $C' = ab$ (with, as always, $a < b$ with respect to the order for $\text{bot}(D)$), which is nonmonotone with respect to the characteristic subdomain $D$, which contains $ef$ as a subarc and for which $b \le b_1$. Clearly, $bf \cap ef = \{f\}$, so that $bf \preceq ef$. If we let
$$c' = \sup\{z : z \in \text{bot}(D) \cap bf, z \le b_1\}$$
and
$$c'' = \inf\{z : z \in \text{bot}(D) \cap bf, z \ge b_2\},$$
then $L = c'c'' \in \mathcal{I}(D)$ so that from the minimality of $ef$ it follows that $c' \ge a_2$ and we have conclusion (*iii*). In exactly the same manner one obtains conclusion (*ii*) in the case that $s \in (0, \lambda_0)$, $l(s) > a_1$ for all $s \in (0, \lambda_0)$.∎

**Comment 3.18.** It is to be noted that Case (*ii*) of the conclusion of the preceding proposition does not exclude the possibility that the nonmonotone extended $i$-characteristic $ab$ also contains a subarc $L' \in \mathcal{I}(D)$ for which $B \preceq L' \preceq U$, that is, that Cases (*ii*) and (*iii*) occur simultaneously. Nor does the conclusion preclude $ab$ having many proper contact points on either or both sides of $A$. Similar comments apply to Case (*iii*). In Case (*i*), $ef$ has to have at least one proper contact point between $A$ and $B$, but it could have many such points as well as proper contact points to the left of $A$ or the right of $B$.

**Proposition 3.19.** (Limit Characteristic Principle) *Let $(D, C_0)$ be an $i$-characteristic subdomain. Let $p \in \text{bot}(D)$. If there is no nontrivial exiting extended $i$-characteristic at $p$, then there is a sequence $\{C_k\}$ of elementary $i$-characteristics for which*
(*i*) $\{p\} \preceq C_k$
(*ii*) $\text{diam}(C_k) \to 0$.

**Proof**. Assume that there is no nontrivial exiting extended $i$-characteristic at $p \in \text{bot}(D)$. Let $E = ef = \min_D(\{p\})$. Assume that $e \neq f$. Obviously, $p \notin E$, so that $E \cap \text{bot}(D) = \{e, f\}$, since otherwise $E$ would have a proper subarc $E' \in \mathcal{C}(\{p\})$, which would contradict the fact that $E = \min_D(\{p\})$. Thus $E$ is an elementary $i$-characteristic. We now proceed as in the preceding proof. Let $w_0$ be any point of $E$ other than $e$ or $f$, and consider a small $j$-arc $J \subset G$, which extends from $w_0$ to the right of $E$ (i.e., into $D$). Let $J$ be parametrized by $z = w(s)$, $0 \leq s \leq \lambda_0$ with $w(\lambda_0) = w_0$, and let $C(s) = l(s)r(s)$ be the elementary $i$-characteristic through $w(s)$. As in the preceding proof, it follows from the minimality of $E$ that either $l(s) > p$ for all $s \in (0, \lambda_0)$ or $r(s) < p$ for all $s \in (0, \lambda_0)$. Assume, for definiteness, that the latter occurs. By considering what happens as $s \to 0$ we obtain, as in the preceding proof, a nonmonotone extended characteristic $ab$ (with $a < b$ with respect to the order on $\text{bot}(D)$) which contains $ef$ as a proper subarc, for which $b < p$ and for which $bf \preceq ef$. But $p \notin bf$ by hypothesis, so that $bf$ must have a subarc $gh$ which is an elementary $i$-characteristic for which $g < p < h$, which contradicts the minimality of $E$. This contradiction implies that we must in fact have $e = f$, which in turn implies, in light of the first conclusion of Proposition 3.14 that the desired sequence $\{C_k\}$ exists. ∎

**Definition 3.20.** A point $p \in \partial G$ is called a *regular point* of $\theta$ if $\theta(z)$ has a nontangential limit as $z \to p$; otherwise it is called a *singular point*.

**Proposition 3.21.** *If $p \in \partial G$ is a regular boundary point, then there must be at least one elementary characteristic which exits at $p$ and which makes an angle $\phi \in (0, \pi)$ with $\partial G$ at $p$.*

**Proof**. This follows immediately from the Peano existence for (local) solutions of the initial value problem $y' = F(x, y)$, $y(x_0) = y_0$ when $F$ is continuous in a neighborhood of $(x_0, y_0)$.∎

**Definition 3.22.** Let $p \in \partial G$ be a singular boundary point of $\theta$, and let $z = z(s)$, $0 \leq s \leq l$ be the parametrization of an nontrivial arc of an extended characteristic with $z(0) = p$. If $\lim_{s \to 0} \theta(z(s))$ exists, then $\theta$ is said to have a *singularity of type 0* at $p$. Otherwise $\theta$ is said to have a *singularity of type 1* at $p$.

**Definition 3.23.** If $\theta$ has no nontrivial exiting extended characteristic at $p$ then $\theta$ is said to have a *singularity of type 2* at $p$.

The following proposition is needed to make use of fans of characteristics in what is to follow.

**Proposition 3.24.** *Let $C_1 = pq_1$ and $C_2 = pq_2$ be distinct nontrivial elementary $i$-characteristics. Let $A$ be the arc of $\partial G$ joining $q_1$ and $q_2$ which does not contain $p$ and let $T$ denote the interior of the curvilinear triangle bounded by $C_1 \cup C_2 \cup A$. Then one of the following must happen*

*(i) There is some $d > 0$ such that for all $z \in N(p,d) \cap T$ the elementary $i$-characteristic through $z$ exits at $p$.*

*(ii) $p$ is proper contact point of some extended characteristic $C$.*

**Proof.** Let $W$ be the set of all $z \in T$ for which the elementary $i$-characteristic $C(z)$ through $z$ doesn't exit at $p$. For $z \in W$, $C(z)$ has no points in common with either $C_1 \setminus \{q_1\}$ or $C_2 \setminus \{q_2\}$, so that both endpoints of $C(z)$ lie on $A$. For each $z \in W$, let $D(z)$ denote the interior of the domain bounded by $C(z)$ and the arc of $A$ whose endpoints are those of $C(z)$. For $w_1, w_2 \in T$, $C(w_1) \cap T$ and $C(w_2) \cap T$ are either identical or disjoint, so that for $w_1, w_2 \in W$, $D(w_1)$ and $D(w_2)$ are either nested or disjoint. Let $\xi = \frac{1}{2}\mathrm{dist}(p,A) > 0$. If $z \in N(p, \frac{\xi}{2}) \cap W$, then it follows from Proposition 2.8 that $\mu(D(z)) \geq \eta \xi^2$ since $\lambda(C(z) \cap N(p,\xi)) \geq \xi$ and each of the $j$-arcs joining a point of $C(z) \cap N(p,\xi)$ to $A$ in $D(z)$ has length at least $\xi$. Thus if $z_1, \ldots z_n \in N(p, \frac{\xi}{2}) \cap W$ are such that the corresponding $D(z_k)$ are disjoint, then $n \leq \frac{\mu(G)}{\eta \xi^2}$. If $(i)$ is not true, a pigeonhole (area exhaustion) argument then shows that there is a sequence $\{z_k : k \geq 1\}$ of points of $W$ tending to $p$ such that $\{D(z_k)\}$ is an increasing sequence. But our extended characteristic construction then gives us a monotone extended characteristic $C$ (with respect to the characteristic subdomain bounded by $C(z_1)$ and the complement in $\partial G$ of the subarc of $A$ joining the endpoints of $C(z_1)$) for which $p \in C$. But the endpoints of $C$ are in $A$ since those of all the $D(z_k)$ are in $A$, so that $p$ is a proper contact point of $C$. Thus $C$ satisfies $(ii)$. ∎

**Definition 3.25.** ($i$-fan) A solution $\theta$ is said to have an $i$-fan at $p \in \partial G$ if $p$ is not a proper contact point of an extended characteristic and there is more than one elementary $i$-characteristic of $\theta$ exiting at $p$. In this case the family of all elementary $i$-characteristics exiting at $p$ is denoted by $\mathcal{F}_i(p)$. The point $p$ is referred to as the *vertex* of the fan.

Proposition 3.24 implies that if $C_1$ and $C_2$ are distinct elementary $i$-characteristics in $\mathcal{F}_i(p)$, then there is a $d > 0$ such that all points of $N(p,d)$ between them belong to members of $\mathcal{F}_i(p)$, so that the interior of $\cup \{C : C \in \mathcal{F}_i(p)\}$ is nonempty. Furthermore, it follows from Proposition 3.7 that the $j$-characteristic through any point in the interior of $\cup \{C : C \in \mathcal{F}_i(p)\}$ is strictly concave towards the side facing $p$. These two facts in turn imply that if there are $i$-fans at $p_1 \neq p_2$ then the interiors of $\cup \{C : C \in \mathcal{F}_i(p_k)\}$, $k = 1, 2$ must be disjoint. Since $\mu(\cup \{C : C \in \mathcal{F}_i(p)\}) > 0$, *the set of points which are the vertices of $i$-fans is at most countable*.

**Proposition 3.26.** *Let there be a regularly exiting extended characteristic and no fan at $p$. Then $p$ is a regular point. Moreover, if $\theta$ is the inclination function of a locally Lipschitz solution in $G$ of a normal system, then the corresponding $R = (R_1, R_2)$ has a nontangential limit at $p$.*

**Proof.** Let $F$ be an extended $i$-characteristic which exits regularly at $p$. Without loss of generality we can assume that the positive direction along $\partial G$ at $p$ is that of the positive real axis. Let $C$ be a nontrivial arc of $F$ parametrized by $z = z(s)$, $0 \leq s \leq L$, with $z(0) = p$. Again without loss of generality we can assume that $\arg\{z'(s)\} = \theta(z(s))$. Let $\gamma = \arg\{z'(0)\}$. We consider the following three possibilities separately.

*(i)* $\gamma \in (0, \frac{\pi}{2}) \cup (\frac{\pi}{2}, \pi)$. In this case we can assume that $C$ is an arc of an elementary $i$-characteristic exiting at $p$. To be specific, we assume that $\gamma \in (0, \frac{\pi}{2})$. Let $q$ be an interior point of $C$. For any $\epsilon > 0$ there exist nontrivial $j$-characteristic arcs $E^+ = E^+(\epsilon)$ and $E^- = E^-(\epsilon)$ emanating from $q$ to the right and left of $C$ such that $\lambda(\theta(E^+ \cup E^-)) < \epsilon$. Since there is no fan at $p$, none of the $j$-characteristics through any point of $E^+ \cup E^-$ other than $q$ exits at $p$. Let $C(\epsilon)$ be an initial segment of $C$ such that $\lambda(\theta(C\setminus\{p\})) < \epsilon$. Let $\alpha < \gamma$ we show that $\theta(z) \to \gamma$ as $z \to p$ between $C$ and the ray $\arg\{z\} = \alpha$. It is easy to see that some initial segment of this ray (that is, the portion of the ray contained in $N(p,\xi)$ for some $\xi > 0$) is covered by a collection $\mathcal{Q}(\epsilon)$ of characteristic quadrilaterals $Q$ one of whose $i$-sides is a subarc of $C(\epsilon)$ and one of whose $j$-sides is a translate of an initial subarc of $E^+(\epsilon)$. The quasi-HP property then implies that there is a $\delta = \delta(\epsilon)$ such that if $|z - p| < \delta$ and $z$ is between $C$ and the ray $\arg\{z\} = \alpha$, then $|\theta(z) - \gamma| < 2K\epsilon$. Since $\epsilon$ is arbitrary, $\theta(z) \to \gamma$ as $z \to p$ between $C$ and the ray. It is also clear that for sufficiently small $\epsilon$ the translate of $E^-$ down to $p$ is a nontrivial initial arc of a $j$-characteristic $C'$ which exits regularly at $p$ and which forms with $\partial G$ an acute angle of size $\frac{\pi}{2} - \gamma$. What we have shown in regard to $C$ now implies that if $\beta < \frac{\pi}{2} - \gamma$, then $\theta(z) + \frac{\pi}{2} \to \gamma + \frac{\pi}{2}$ as $z \to p$ between $C'$ and the ray $\arg\{z\} = \pi - \beta$, so that $\theta(z) \to \gamma$ as $z \to p$ in that curvilinear sector. Finally, it is easy to see that $\theta(z) \to \gamma$ as $z \to p$ between $C$ and $C'$. Since $\alpha < \gamma$ and $\beta < \frac{\pi}{2} - \gamma$ are arbitrary, we have the desired regularity at $p$. The case $\gamma \in (\frac{\pi}{2}, \pi)$ is handled in the same manner apart from minor changes of a notational nature. In the case that $\theta$ arises from the solution of a normal system it is clear from (*ii*) of Definition 1.1 that both $R_i(z)$ and $R_j(z)$ have limits as $z \to p$ in the sector $\alpha < \arg\{z\} < \beta$, so that the second conclusion is valid in this case.

(*ii*) $\gamma = \frac{\pi}{2}$. Here again we can take $C$ to be an arc of an elementary $i$-characteristic exiting at $p$, and proceed as in case (*i*). For any $\alpha < \frac{\pi}{2}$ it follows as in case (*i*) that $\theta(z) \to \gamma$ as $z \to p$ between $C$ and the ray $\arg\{z\} = \alpha$. Again it is immediate that the same holds in the curvilinear sector between $C$ and the ray $\arg\{z\} = \pi - \alpha$. The second conclusion is likewise immediate.

(*iii*) $\gamma \in \{0, \pi\}$. We deal with the case $\gamma = 0$, the case $\gamma = \pi$ being essentially the same. Here, we can define $E^-(\epsilon)$ as in case (*i*), but because there is no fan at $p$ the translate of $E^-(\epsilon)$ down to $p$ is a nontrivial arc of a $j$-characteristic which is orthogonal to $\partial G$ at $p$. This puts us in case (*ii*), so we are done. ∎

**Proposition 3.27.** *If the extended $j$-characteristic $C$ exits singularly at $p$, then there is no nontrivial exiting extended $i$-characteristic at $p$.*

**Proof**. Let $C$ be given parametrically by $z = z(s)$, $0 \le s \le L$, with $z(0) = p$, and let $\phi(s) = \arg\{z'(s)\}$. For $0 < s < L$ let $E(s)$ be the elementary $i$-characteristic containing $z(s)$. By Proposition 3.15, $\mathrm{diam}(E(s)) \to 0$ as $s \to 0$. If for some $s_0$ an endpoint of $E(s_0)$ were to coincide with $p$ then the same would be the case for $E(s)$ for all $s \in (0, s_0)$, giving us an $i$-fan at $p$. But this is inconsistent with the fact that the sets $P(s_0, T)$ and $N(s_0, T)$ of Proposition 3.15 have positive measure, since a $j$-characteristic joining points on two distinct elementary $i$-characteristics of an $i$-fan must be strictly concave towards the vertex of the fan. Thus both of the endpoints of all the $E(s)$ are different from $p$. Say now that to the contrary there is a nontrivial exiting extended $i$-characteristic $I$ at $p$, the distance between whose endpoints is $d > 0$. Then from the definition of extended characteristic it would follow that for any $\epsilon > 0$ there is an elementary $i$-characteristic $I'_\epsilon$ the distance between whose endpoints is at least $d/2$, and for which $\mathrm{dist}(p, I'_\epsilon) < \epsilon$. But since $\mathrm{diam}\, E(s) \to 0$, $I'_\epsilon \ne E(s)$ for $s \in (0, s_0)$ for some $s_0 > 0$. But then for any sufficiently small $s$, $I'_\epsilon$ would have to have a point in common with $E(s)$, an obvious contradiction. ∎

**Definition 3.28.** A singular point $p \in \partial G$ at which there is no nontrivial exiting extended $i$-characteristic is called an *$i$-singularity*.

The following is an immediate consequence of Propositions 3.26 and 3.27 and the preceding definition.

**Corollary 3.29.** *If $p \in \partial G$ is a singular point, then at least one of the following must be the case.*
*(i) $p$ is the vertex of a fan*
*(ii) $p$ is an 1-singularity*
*(iii) $p$ is a 2-singularity*

**Comment 3.30.** The proof of the Main Theorem 2.2 in the following section is based on showing that the set of $i$-singularities at which there is no fan has $\alpha$-dimensional Hausdorff measure 0 for some $\alpha < 1$.

**Proposition 3.31.** (Type 1 Singularities with Fans) *Let $C$ be a subarc of an elementary $i$-characteristic which joins $q \in G$ to $p \in \partial G$ and which is parametrized by $z(s)$, $0 \leq s \leq l$ with $z(0) = p$. Assume that for some $\eta < \frac{\pi}{2}$*

(3.22) $$\phi + \eta \leq \arg\{z'(s)\} \leq \phi + \pi - \eta, \text{ for } 0 < s \leq l,$$

*where $\phi$ is the argument of the positive tangent to $\partial G$ at $p$. If $\lim_{s \to 0} \arg\{z'(s)\}$ does not exist, then $C$ belongs to a nontrivial $i$-fan with vertex $p$.*

**Proof.** For $0 < \xi < \pi$, let $R(p, \xi)$ denote the (open) ray
$$R(p, \xi) = \{z : \arg\{z - p\} = \xi\}.$$
Without loss of generality we can assume that $\phi = 0$. It follows immediately from hypothesis (3.22) that $C$ lies in the sector $\{z : \eta \leq \arg\{z - p\} \leq \pi - \eta\}$. The hypothesis that $\lim_{s \to 0} \arg\{z'(s)\}$ doesn't exist implies that there is some $\delta$, $0 < \delta \leq \frac{\eta}{K+1}$, such that for arbitrarily small $\sigma_1 > \sigma_2 > 0$,

(3.23) $$|\arg\{z'(\sigma_1)\} - \arg\{z'(\sigma_2)\}| \geq 2\delta.$$

For $0 < s \leq l$, let $F(s) = a(s)b(s)$ be the elementary $j$-characteristic through $z(s)$. It follows from Proposition 3.27 that $a(l) \neq p \neq b(l)$. Let $D$ be the $j$-characteristic subdomain bounded by $F(l)$ and the arc of $\partial G$ joining $a(l)$ to $b(l)$ and containing $p$. Clearly, $F(s) \cap C \cap G = \{z(s)\}$ and $a(s)$ is nonincreasing and $b(s)$ is nondecreasing with respect to the order on $\text{bot}(D)$. Since by Proposition 3.15, $\text{diam}(F(s)) \to 0$ as $s \to 0$, by replacing $C$ by a sufficiently short initial subarc we can assume that

(3.24) $$R(p, \xi) \cap \partial D \subset F(l), \frac{\delta}{2} \leq \xi \leq \pi - \frac{\delta}{2}.$$

Now consider any pair of numbers $\sigma_1 \neq \sigma_2$ in $(0, l)$ for which (3.23) holds. Obviously, for at least one of them, call it $\sigma$, there must hold $|\arg\{z'(\sigma)\} - \frac{\pi}{2}| \geq \delta$. Then $F(\sigma)$ must have a subarc $W = uv$ such that $\lambda(\theta(W)) = \delta$ and $z(\sigma) \in \{u, v\}$. To see this, say, for example, that $\arg\{z'(\sigma)\} \leq \frac{\pi}{2} - \delta$. If along the arc of $F(\sigma)$ emanating to the left of $C$, the argument of the tangent were always within $\delta$ of its argument at $z(\sigma)$, then $F(\sigma)$ would never cross $R(p, \pi - \delta)$, and would therefore not be able to exit $D$. An analogous argument may be used in the case that $\arg\{z'(\sigma)\} \geq \frac{\pi}{2} + \delta$.

From this it follows that there is a sequence $s_k \to 0$ of such numbers $\sigma$ for all of which the corresponding $j$-arc $W_k = u_k v_k$ lies on one side of $C$ or the other. To be specific, say $u_k = z(s_k)$ and $W_k$ lies to the right of $C$. It follows from the quasi-HP-property that if $W' = z(s)v$ is any translate of $W_k$ along $C$ in the direction of increasing $s$, then
$$\delta/K \leq \lambda(\theta(W')) \leq K\delta.$$
From this and the fact that $\delta \leq \frac{\eta}{K+1}$ it follows in turn that for $s > s_k$, along any $i$-characteristic arc whose initial point is in $W_k$ and which is parallel to the subarc $z(s_k)z(s)$ of $C$, the inclination of the tangent is in the interval
$$[\eta - K\delta, \pi - \eta + K\delta] \subset [(K+1)\delta - K\delta, \pi - (K+1)\delta + K\delta] = [\delta, \pi - \delta].$$
For each $k$ let
$$S_k = \{R(z, \xi) : z \in W_k \text{ and } \xi \in [\delta, \pi - \delta]\}.$$
Now, $\operatorname{diam}(W_k) \leq \operatorname{diam}(F(s_k))$ and the latter tends to 0 as $k \to \infty$, so that in light of (3.24) by eliminating a finite number of elements of $\{s_k\}$ and renaming we can assume that all the $R(z, \xi)$ in each of the $S_k$ meet $\partial D$ a point of $F$ before exiting $G$. From this it follows that the translates of $z(s_k)z(l)$ along $W_k$ all lie in $G$. This means that for any $s \in [s_k, l]$ the complete translate of $W_k$ along $C$ from $z(s_k)$ up to $z(s)$ belongs to $G$. For the translate $T_k$ of $W_k$ up to $z(l)$ we have that $z(l) \in T_k \subset F(l)$, and by the quasi-HP property $\lambda(\theta(T_k)) \geq \delta/K$. This bound implies that there is a positive lower bound on the length of the $T_k$. Thus $\cap \{T_k : k \geq 1\} = T$, is a nontrivial arc of $F(l)$ one of whose endpoints is $z(l)$. However, the translate of $T$ down to $z(s_k)$ is contained in $W_k \subset F(s_k)$. Since $\operatorname{dist}(F(s), \{p\}) \to 0$ as $s \to 0$, it follows that all $i$-characteristics through $T$ exit at $p$, and this establishes the existence of the desired fan.∎

We need the following elementary

**Lemma 3.32.** *Let $f$ be a bounded measurable function on $[0, T]$ such that $\int_0^T f(x)dx = A > 0$. Then for any $\rho \in [0, \frac{1}{2})$*
$$\lambda(\{\xi : f(\xi) \geq \rho \tfrac{A}{T} \text{ and } \int_0^\xi f(x)\,dx \geq \rho A\}) > 0.$$

**Proof.** Let $f(x) \leq M$ a.e. on $[0, T]$ and let $Q = \{\xi : f(\xi) \geq \rho \tfrac{A}{T}\}$. Obviously, $\lambda(Q) > 0$. Let $Q^*$ be the set of density points of $Q$, so that $\lambda(Q^*) = \lambda(Q) > 0$. Let $\xi_0 = \sup\{\xi : \xi \in Q^*\}$. If the set $Q^* \cap \{\xi : \int_0^\xi f(x)\,dx \geq \rho A\}$ had measure 0, then for each $\delta > 0$ there would be a point $\xi \in (\xi_0 - \delta, \xi_0)$ for which $\int_0^\xi f(x)\,dx < \rho A$, so that
$$A = \int_0^\xi f(x)\,dx + \int_\xi^{\xi_0} f(x)\,dx + \int_{\xi_0}^T f(x)\,dx < \rho A + \delta M + (T - \xi_0)\tfrac{\rho A}{T}$$
$$< 2\rho A + \delta M < A,$$
for $\delta$ sufficiently small, so that $\lambda(Q^* \cap \{\xi : \int_0^\xi f(x)\,dx \geq \rho A\})$ must be positive.∎

The following proposition plays a fundamental role in the proof the Main Theorem 2.2.

**Proposition 3.33.** (Essentially Singularity-Free Boundary Arcs) *There is a function $\overline{\eta} = \overline{\eta}(K, \rho) < \rho$, with the following property. Let $G \in \mathcal{G}(\rho)$ and let $\theta \in \operatorname{HP}(G, K)$. Let $C$ be an extended $i$-characteristic of $\theta$ with proper contact point $p \in \partial G$. Let $E = pp'$ be the elementary $j$-characteristic exiting at $p$, and let $D$ be the $j$-characteristic subdomain of $G$ bounded by $E$ and an arc $A$ of $\partial G$. Let $b \in A \cap C$ be the endpoint of $C$ on $A$ and let $F$ be the subarc $[b, p']$ of $A$. If $\operatorname{dist}(p, F) \geq \rho$, then the subarc $[p, p'']$ of $A$ with $|p'' - p| = \overline{\eta}$ has at most countably many boundary singularities of $\theta$.*

**Proof.** Without loss of generality we may assume that $p = 0$, that the positively oriented tangent to $\partial G$ at $p$ has the direction of the positive real axis and that $A$ extends to the right of $p$. Let $C$ be parametrized by $z(s) = x(s) + iy(s)$, with $z(0) = p$. We shall introduce positive constants, which will depend only on $K$ and $\rho$. In most cases these constants will be denoted by the same symbol $B$, but the use of a single symbol to denote different constants will not cause confusion or be misleading since any statement in which $B$ plays a role will be valid if $B$ is taken to be any suitably large number and since this convention is used only finitely many times. The symbol $B_1$, on the other hand, will refer to a specific constant which again depends only on $K$ and $\rho$, and will have the same value every time it is used.

If we give $A$ in nonparametric form by $y = g(x)$, then
$$(3.25) \qquad |g(x)| \leq \tfrac{x^2}{\rho} \text{ and } |g'(x)| \leq \tfrac{2x}{\rho}, \text{ for } |x| \leq \tfrac{1}{10\rho}.$$
Clearly, Proposition 3.12 holds for extended characteristics with the obvious wording changes. Any $j$-arc emanating to the left of $C$ from any point $z$ of the subarc $pb$ of $C$ will exit $G$ at a point of $F$, so that the length of any such $j$-arc is a least $\rho - |p - z|$. Since $\operatorname{dist}(z(s), \partial G) \leq |z(s)| \leq |s|$, it follows from Proposition 3.12 that $|\arg\{z'(s)\}| \leq \tfrac{1}{2}$ for $|s| \leq \tfrac{1}{B}$ so that for these values of $s$, we have
$$\tfrac{1}{B} \leq x'(s) \leq 1 \text{ and } |y'(s)| \leq 1,$$
and can represent $C$ nonparametrically by $y = f(x)$, where $f(x(s)) = y(s)$. We have
$$|f'(x(s))| = |\tan(\arg\{z'(s)\})| \leq 2|\arg\{z'(s)\}| \leq B\sqrt{\operatorname{dist}(z(s), A)}$$
$$\leq B\sqrt{\operatorname{dist}(z(s), \mathbb{R}) + \operatorname{dist}(x(s), A)} \leq B\sqrt{y(x(s)) + x(s)^2},$$
where we have used the first bound of (3.25) in the form $\operatorname{dist}(x(s), A)) \leq \tfrac{(x(s))^2}{\rho}$. Thus there is a $B_1 \geq 1$ such that
$$(3.26) \qquad |f'(x)| \leq B_1\sqrt{|f(x)| + x^2} \text{ and } |f(x)| \leq B_1|x|, \text{ for } |x| \leq \tfrac{1}{B_1}.$$
For $|x| \leq \tfrac{1}{B_1}$ we then have the following. In the first place,
$$|f'(x)| \leq B_1\sqrt{B_1|x| + x^2} \leq 2^{1/2}B_1^{3/2}|x|^{1/2},$$
so that in fact $|f(x)| \leq \tfrac{2}{3}2^{1/2}B_1^{3/2}|x|^{\tfrac{3}{2}} < B_1^{3/2}|x|^{\tfrac{3}{2}}$. Repeating this argument we have that $|f'(x)| \leq 2^{1/2}B_1^{7/4}|x|^{\tfrac{3}{4}}$, so that in fact $|f(x)| \leq \tfrac{4}{7}2^{1/2}B_1^{\tfrac{7}{4}}|x|^{\tfrac{7}{4}} < B_1^{\tfrac{7}{4}}|x|^{\tfrac{7}{4}}$, so that $|f'(x)| \leq 2^{1/2}B_1^{\tfrac{15}{8}}|x|^{\tfrac{7}{8}}$, and so on. Thus, in fact
$$(3.27) \qquad |f(x)| \leq B_1^2 x^2, \text{ for } |x| \leq \tfrac{1}{B_1},$$
and
$$(3.28) \qquad |f'(x)| \leq \sqrt{2}B_1^2|x|, \text{ for } |x| \leq \tfrac{1}{B_1}.$$

Before proceeding we remind the reader that each time the symbol $B$ appears, it may have a value larger than it did at its previous appearence. For any $t \in [0, \tfrac{1}{2B}]$, let $J(t)$ be the $j$-half-characteristic emanating to the right of $C$ (that is, downwards) from $t + if(t) = \alpha(t)$ and joining it to a point of $A$. Note that $J(t)$ will reduce to a single point if $\alpha(t) \in A$. Let $J(t)$ be given parametrically by $\zeta_t(s)$ with $\zeta_t(0) = \alpha(t)$. Let $0 < x_0 \leq \tfrac{1}{B}$ and let $t_0 = x_0^2$. Then $|f(t_0) - g(t_0)| \leq Bx_0^4$. For notational convenience we can assume without loss of generality that $\arg\{z'(s)\} = \theta(z(s))$. We show that, with an appropriately large value of $B$, $0 < x_0 \leq \tfrac{1}{B}$ implies
$$(3.29) \qquad |\theta(\zeta_{t_0}(s_1)) - \theta(\zeta_{t_0}(0))| \leq x_0 = \sqrt{t_0},$$
for every $\zeta_{t_0}(s_1) \in J(t_0)$. If this were not true, there would be an $s_1 \in (0, \lambda(J(t_0)))$ such that $|\theta(\zeta_{t_0}(s_1)) - \theta(\zeta_{t_0}(0))| = x_0$. Let $s_2$ be the smallest number in $[0, s_1]$ for which $|\theta(\zeta_{t_0}(s_2)) - \theta(\zeta_{t_0}(0))| = x_0$, so that

(3.30) $$\lambda(\theta(\zeta_{t_0}([0, s_2]))) \leq 2x_0.$$
By increasing $B$ if necessary $t_0 = x_0^2 < \frac{1}{B^2}$ will be so small that
(3.31) $$s_2 \leq 2|f(t_0) - g(t_0)| \leq Bx_0^4.$$
By the preceding lemma (with $\rho = \frac{1}{3}$, $T = s_2$, $A = x_0$) it then follows that there is an $s_3 \in (0, s_2)$ for which
$$\left|\frac{d\theta(\zeta_{t_0}(s))}{ds}\right|_{s=s_3}\right| \geq \frac{x_0}{3s_2} \geq \frac{1}{Bx_0^3}$$
and
(3.32) $$|\theta(\zeta_{t_0}(s_3)) - \theta(\zeta_{t_0}(0))| \geq \frac{x_0}{3},$$
where the expressions inside the absolute values on the left side of these last two bounds have the same sign. Let $R$ be the $i$-half-characteristic through $\zeta_{t_0}(s_3)$ emanating from $J(t)$ towards its concave side (that is, to the left of $J(t_0)$ if $d\arg\{\zeta'_{t_0}(s)\}/ds$ is positive at $s_3$, and to the right if it is negative). It follows from the characteristic length bound (Proposition 2.5) that
(3.33) $$\lambda(R) \leq Bx_0^3.$$
We concentrate on the case that that $R$ emanates to the left of $J(t_0)$; the opposite case is much easier to handle as we indicate below. By (3.27), (3.31) and (3.33), $R$ exits $G$ at a point $g(\beta)$ with
(3.34) $$\beta \leq t_0 + |f(t_0)| + |\zeta_{t_0}(0) - \zeta_{t_0}(s_3)| + \lambda(R)$$
$$\leq x_0^2 + O(x_0^3) \leq 2x_0^2, \text{ for } x_0 \leq \frac{1}{B}.$$
For $t \geq t_0$ let $\overline{J}(t)$ be the $j$-arc joining $\alpha(t) \in C$ to a point $r(t) \in R$, so that $\overline{J}(t_0) = \zeta_{t_0}([0, s_3])$. Obviously, $\overline{J}(t)$ is not defined for all $t \in [t_0, \frac{1}{2B_1}]$, but for any $t$ for which it is defined, $\overline{J}(t_0)$ and $\overline{J}(t)$ are the $j$-sides of a characteristic quadrilateral $Q(t)$, whose $i$-sides are the arc $\alpha(t_0)\alpha(t)$ of $C$ and the arc $\zeta_{t_0}(s_3)r(t)$ of $R$. From the $K$-quasi-HP property and (3.30) it follows that
(3.35) $$\lambda(\theta(\overline{J}(t))) \leq K\lambda(\theta(\overline{J}(t_0))) \leq K\lambda(\theta(\zeta_{t_0}([0, s_2]))) \leq 2Kx_0,$$
and by (3.28)
(3.36) $$\lambda(\theta(\zeta_{t_0}(s_3)r(t))) \leq K\lambda(\theta(\alpha(t_0)\alpha(t))) \leq 2K\tan^{-1}(Bt) \leq Bt.$$
But for any $t \leq 10x_0^2$ for which $J(t)$ is defined, we have from (3.35) and (3.28) that
$$|\theta(\zeta_t(s))| \leq \theta(\alpha(t)) + 2Kx_0 \leq \tan^{-1}(Bt) + 2Kx_0$$
$$\leq 20Bx_0^2 + 2Kx_0 \leq \frac{1}{100},$$
for $x_0 \leq \frac{1}{B}$. From this together with (3.34) it is clear that $g(\beta)$ is the endpoint of $\overline{J}(t_1)$ for some $t_1 \in [x_0^2, 10x_0^2]$. It follows from the quasi-HP-property that $R$, when oriented from $\zeta_{t_0}(s_3)$ to $r(t_1)$ has a well defined (one-sided) tangent whose argument we denote by $\phi_0$. But then we have from (3.28), (3.32), (3.36) and the quasi-HP-property that
$$\phi_0 \geq \frac{x_0}{3} - |\theta(\alpha(t_0))| - Bt_1 \geq \frac{x_0}{2} - \theta(\alpha(x_0^2)) - 10Bx_0^2 \geq \frac{x_0}{3} - Bx_0^2,$$
so that
(3.37) $$\phi_0 \geq \frac{x_0}{4}.$$
for $x_0 \leq \frac{1}{B}$. But since by (3.25) and (3.34)
$$|\tan^{-1}(g'(\beta))| \leq |g'(\beta)| \leq \frac{2\beta}{\rho} \leq \frac{4}{\rho}x_0^2, \text{ for } x_0 \leq \frac{1}{B}.$$
Taking intt account the direction from which $R$ crosses $A$ at $g(\beta)$, one sees that that $\phi_0 \leq |\tan^{-1}(g'(\beta))|$, which is a contradiction for $x_0 < \min\{\frac{1}{B}, \frac{\rho}{15}\}$. If $R$ emanates to the right of $J(t_0)$ we get the same contradiction more easily because in that case we move

back along $C$ towards $p$ and there are smaller bounds throughout. This establishes (3.29) for $x_0 \leq \frac{1}{B}$ (with, of course, an appropriately large $B$).

Thus for $\sqrt{t} \leq \frac{1}{B}$, $\lambda(\theta(J(t))) \leq \sqrt{t}$. From this together with (3.27) it follows that

(3.38) $\qquad |\theta(z)| \leq \sqrt{t} + B|t| \leq \frac{1}{100}$, for $z \in J_t$, $0 < t < \frac{1}{B}$.

Finally, after increasing $B$ if necessary, we have $|\tan^{-1}(g'(t))| \leq \frac{1}{100}$ for $|t| \leq \frac{1}{B}$. Let $l_t = \lambda(J_t)$ and let $q(t) = \zeta_t(l_t) \in A$. We now show that $\{q(t) : 0 \leq t \leq \frac{1}{B}\}$ is an arc of $A$, for which we only have to show that $q(t)$ is continuous on $(0, \frac{1}{B})$. Clearly we have from (3.38) that

(3.39) $\qquad |\arg\{\frac{\partial \zeta_t(s)}{\partial s}\} + \frac{\pi}{2}| \leq \frac{1}{100}$, for $t \leq \frac{1}{B}$ and $0 \leq s \leq \lambda(t)$.

If $l_{t_0} = 0$, then what we have already shown implies that $l_t \to 0$ as $t \to t_0$, from which it immediately follows that $q(t) \to q(t_0)$ as $t \to t_0$. Thus consider a fixed $t_0$ for which $l_{t_0} > 0$. Let $\delta > 0$ be small and let $s \in (l_{t_0} - \delta, l_{t_0})$. Let the elementary $i$-characteristic through $\zeta_{t_0}(s)$ parametrized by $z = w_s(\sigma)$ with $w_s(0) = \zeta_{t_0}(s)$. It is clear that there is some $\delta' \in (0, \delta]$ such that for all $\sigma \in [-\delta', \delta']$, $w_s(\sigma) \in J(t(\sigma))$ (where $t(\sigma)$ is some continuous function of $\sigma$) and $|\theta(w_s(\sigma)) - \theta(w_s(0))| \leq \frac{1}{100}$. Then it is clear from (3.39) that for $|\sigma| \leq \delta'$, $J(t(\sigma))$ intersects $A$ at a point within $2\delta$ of $q(t_0)$. The desired continuity follows immediately since there is a positive $\delta''$ such that for $|t - t_0| < \delta''$, $J(t)$ intersects $w_s((-\delta', \delta'))$.

From this and (3.38) it follows that in fact for $0 < t < \frac{1}{B}$, $|\theta(z)| \leq 2\sqrt{t}$ for $z$ in the domain bounded by $J_t$, and the arcs of $C$ and $A$ with endpoints $p$ and $q(t)$. The desired conclusion now follows immediately from Propositions 3.26 and 3.31 and the fact that at most a countable set of points on $\partial G$ can be the vertices of fans of $\theta$. ∎

We end this section with a technical proposition to be used in the proof of the Main Theorem given in Section 4. The numbers $\eta = \eta(K)$ and $\eta' = \eta'(K)$ arising here are those of Propositions 2.8 and 2.9, respectively.

**Proposition 3.34.** (Trapped Area Bound) *There is a positive number $\eta_1 = \eta_1(K) \leq \min\{1, \eta(K)\}$ with the following property. Let $(D, V)$ be an $i$-characteristic subdomain. Let $U = u_1 u_2 \in \mathcal{I}(D)$ be an elementary $i$-characteristic such that $|p - q| \geq \xi$ for any $p, q \in \text{bot}(D)$ for which $[u_1, u_2] \subset [p, q]$. Assume furthermore that $U$ is a subarc of an extended characteristic lying in $\overline{D}$ of which $u \in \{u_1, u_2\}$ is a proper contact point. Let $U' \subset \overline{D}$ be an elementary $j$-characteristic lying between $U$ and $V$ whose endpoints are $u$, $u'$, with $u' < u_1$ or $u' > u_2$ according as $u = u_1$ or $u = u_2$, respectively (in other words, $u$ lies between $u'$ and the other endpoint of $U$). Assume, furthermore, that $|\Delta\theta(U')| \geq \frac{\pi}{4}$. Then the area of the subdomain of $D$ between $U$ and $V$ is bounded below by $\eta_1 \xi^2$.*

**Proof**. For definiteness we assume that $u = u_1$, so that $u' < u_1$. Let $V = v_1 v_2$. Let $U'$ be parametrized by $w = w(s)$, $0 \leq s \leq L$ with $w(0) = u_1$, and for $s \in [0, L)$ denote by $U(s) \preceq V$ the elementary $i$-arc containing $w(s)$. One of the endpoints $a(s)$ of $U(s)$ lies below $U'$ and the other endpoint $b(s)$ is in $[v_1, v_2] \setminus (u_1, u_2)$. In particular $U(0) = U$. There are the following two possibilities.

Case A. The point $b(s) \geq u_2$ (i.e., lies to the right of $u_2$) for $0 \leq s < L$. In this case $U \preceq U(s)$, and therefore $\lambda(U(s)) \geq \xi$ for all $s \in [0, L)$, so that by Area Bound 2 (Proposition 2.9), the area of $\cup \{U(s) : 0 < s < L\}$ is at least $\eta' \frac{\pi}{4} \xi^2$.

Case B. There is some $s \in (0, L)$ for which $b(s) \leq u_1$. Since $b(0) = u_2$, there is a decreasing sequence $\{s_k\}$ tending to some $\sigma \in [0, L)$ such that $b(\sigma) \geq u_2$, but $b(s_k) \leq u_1$ for $k \geq 1$. But then there is a nonmonotone (with respect to $D$) extended $i$-characteristic $C$ joining $a_0 = \lim_{k \to \infty} a(s_k)$ to a point $b_0 = \lim_{k \to \infty} b(s_k) \leq a_0 \leq u_1$ for which $b(\sigma)$ is a proper contact point. Then $C$ has subarcs $U(\sigma) = a(\sigma)b(\sigma)$ and $W = de$, for which $U(\sigma) \preceq W$ and $U(\sigma) \cap W \subset \{b(\sigma)\}$, and where $d \in [v_1, u']$, $e \in [b(\sigma), v_2]$. Since we therefore have $[u_1, u_2] \subset [d, e]$, we conclude that $\lambda(W) \geq |e - d| \geq \xi$. Regard $W$ as oriented from $e$ to $d$ and consider the subarc $W'$ of $W$ of length $\frac{\xi}{2}$ which starts at $e$, and for each $p \in W'$ consider the $j$-arc $J(p)$ emanating to the left of $W$ and which joins $p$ to a point $q(p) \in \text{bot}(D)$. Because $J(p)$ cannot intersect $C$ at any interior point of $C$ other than $p$, we have $d \leq q(p) \leq u_1$. Since the length of any curve joining $e$ to a point to the left of $u_1$ must be at least $\xi$, it follows that $\lambda(J(p)) \geq \frac{\xi}{2}$. Area Bound 1 (Proposition 2.8) then implies that

$$\mu(\cup \{J(p) : p \in W'\}) \geq \eta \frac{\xi^2}{4}.$$

The conclusions of Cases A and B together imply that the area between $U$ and $V$ is at least $\min\{\eta \frac{\xi^2}{4}, \eta' \frac{\pi}{4} \xi^2\} = \eta_1 \xi^2$, where $\eta_1 = \min\{\frac{\eta}{4}, \frac{\pi}{4} \eta'\}$. ∎

## 4. Proof of the Main Theorem 2.2

The main ingredients in the proof are Proposition 3.33 about essentially singularity free boundary arcs, the limit of characteristic principle (Proposition 3.19), Proposition 3.17 regarding structure of $\min_D(A \cup B)$ and the trapped area bound (Proposition 3.34). We also use the preservation of quasi-HP conditions under linear changes of variable, which allows one to normalize the arc one is working with to have any convenient length. In reference to a quasi-HP function $\theta$ on $G \in \mathcal{G}$, we will say that an arc $E \subset \partial G$ is *essentially $i$-singularity free* (abbreviated ESF-$i$) if it has at most countably many $i$-singularities of $\theta$. As is explained in detail in the final paragraph of the proof, it is enough to prove that there is a number $\tau = \tau(K) < 1$ such that on any suitably small boundary arc the set of all $i$-singularities has Hausdorff dimension at most $\tau$. In the treatment of such arcs there is considerable freedom in the choice of the various explicitly given numerical constants that we use, and for the most part they are far bigger (or smaller) than necessary and have most often been chosen either for the sake of convenience or to avoid the necessity of going into careful geometric arguments. In the same vein, many of the bounds we give are far from being sharp, so that, if we state that the value of some expression is bounded below by 1, it may be clear that it is in fact bounded below by a considerably larger number. Furthermore, in some instances we use symbols to denote numerical constants which could without much effort be determined explicitly.

Let $G \in \mathcal{G}(\rho)$, so that the unsigned curvature of $\partial G$ is everywhere at most $\frac{1}{\rho}$. For $p \in \partial G$ let $w(p) = e^{i\phi(p)}$ be the unit tangent to $\partial G$ (with positive orientation) at $p$. Let $\delta = \delta(G)$ be any positive number small enough that

(1) $\delta < \frac{\rho}{10000}$, which implies in addition that

(2) On any arc $B$ of $\partial G$ of length $100\delta$, $\lambda(w(B)) < \frac{1}{100}$, and that

(3) For any point $p \in \partial G$ and any $r \leq 20\delta$, $\partial N(p,r)$ meets $\partial G$ in two points, so that for such $r$, $N(p,r) \cap G$ is essentially a semidisk.

The corners of $\partial N(p, 9\delta) \cap G$ can be replaced by circular arcs of radius $\frac{\delta}{90}$ and then smoothed at the joins to form a $C^2$ curve which is the boundary of a subdomain
$$\Sigma = \Sigma(p, \delta) \subset N(p, 10\delta) \cap G$$
of $G$ with the following properties:

(1a) $\Sigma \in \mathcal{G}(\frac{\delta}{100})$.

(2a) The arc $N(p, 4\delta) \cap \partial \Sigma$, to be referred to as the *bottom* of $\Sigma$, and which is an arc of $\partial G$, has unsigned curvature bounded above by $\frac{1}{10000\delta}$.

(3a) $\mu(\Sigma) \leq 200\delta^2$.

Now consider the $X = X(p) = \frac{1}{\delta}(\Sigma(p,\delta) - p)$ which has the following properties:

(1b) $X \in \mathcal{G}(\frac{1}{100})$.

(2b) $B = N(p, 4) \cap \partial X$, to be referred as the *bottom* of $X$, is a $C^2$ arc on which the unsigned curvature bounded above by $\frac{1}{10000}$.

(3b) $\mu(X) \leq 200$.

(4b) $\text{dist}(a, \partial X \setminus (N(a,d) \cap B)) \geq \min\{d, \frac{1}{2}\}$ for all $a \in B \cap N(p, 6)$.

Note that what we call the bottom $B$ of $X$ is very close to a straight line segment of length 8 centered at the "midpoint" $p$ of what is close to being the straight portion of the boundary of a semidisk of radius 9. Note also that by Proposition 2.3, $\theta(p + \delta z)$ is a $K$-quasi-HP function on this $G(p)$. A domain $X$ having these four properties will be called a *normalized domain*.

Until the end of the proof of the Main Theorem we work exclusively with quasi-HP functions $\theta$ defined on a normalized domain $X$. It is clear from property (1b) that there is a universal constant $\gamma_0 \in (0,1)$ such that for any such $X$ we have

(4.1) $\qquad |z_1 - z_2| \geq \gamma_0 \text{dist}_{\text{arc}}(z_1, z_2), \quad z_1, z_2 \in \partial X,$

where $\text{dist}_{\text{arc}}(z_1, z_2)$ denotes the length of the shorter arc of $\partial X$ with endpoints $z_1$, $z_2$. We shall prove the following proposition from which, as we will subsequently show, the Main Theorem follows almost immediately. We stress that for an arc $F$ of the bottom of $X$, $\text{diam}(F)$ is the distance between its endpoints.

**Proposition 4.1.** *There is a constant $\Delta_0 = \Delta_0(K) > 0$ such if $\theta$ is a $K$-quasi-HP function on a normalized domain $X$, then the bottom of $X$ has an ESF-i arc of diameter at least $\Delta_0$.*

**Proof**. It is clear from Proposition 2.3 that we can assume that $p = 0$ and that $w(p) = 1$. We establish the desired conclusion as follows (once again, the constants chosen are unnecessarily big or small). Let $H_0$ be an integer for which
$$H_0 \geq \max\{\tfrac{10^6}{\gamma_0^2 \eta_1}, 10000\}$$
and let
$$\epsilon_0 = \overline{\eta}(K, \tfrac{\gamma_0}{2H_0}),$$
where $\eta_1 = \eta_1(K)$ is the constant of the trapped area bound (Proposition 3.34) and $\overline{\eta}(K, \frac{\gamma_0}{2H_0})$ is the constant from Proposition 3.33 about essentially singularity free boundary arcs; all three of these numbers depend solely on $K$. In particular we have

$$\epsilon_0 \leq \tfrac{\gamma_0}{2H_0} < \tfrac{1}{2H_0} \leq \tfrac{1}{20000}.$$

Our strategy will be to show that there is some constant $\gamma_1 = \gamma_1(K)$ such that the assumption that

(4.2) *the bottom $B$ of $X$ has no ESF-$i$ arc of diameter $\epsilon_0$*

leads to the conclusion that the bottom of $X$ has an ESF-$i$ arc of diameter $\tfrac{\epsilon_0}{10000}$ or one of diameter $\gamma_1$, so that $B$ has an ESF-$i$ arc of diameter $\Delta_0 = \min\{\tfrac{\epsilon_0}{10000}, \gamma_1\}$. This underlying assumption (4.2) will be in force throughout and will be used many times.

The proof the existence of such a $\Delta_0$ is divided into three steps.

**Step 1**. There is an elementary $i$-characteristic $C_0 = \overline{a}\overline{b}$ which has the following properties.

(*i*) $|\overline{b} - \overline{a}| \geq 1$.

(*ii*) At least one of $\overline{a}, \overline{b}$ lies in $B$.

**Comment**. We emphasize the words "at least one" since were we able to say "both endpoints lie in $B$", we would be much closer to finishing the proof of the theorem; in effect we could go directly to Step 3. After completing Step 1 we will use what it says to show that there must be a subarc $J$ of $B$, with diam($J$) bounded below by some universal constant $c$ and whose endpoints are joined by an elementary $i$-characteristic; we accomplish this in Step 2 with $c = \tfrac{1}{40}$.

We begin the proof of Step 1 by noting that by (4.2) there is an $i$-singularity within $\epsilon_0$ of the left endpoint of $B$ so that by the limit characteristic principle (Proposition 3.19) there is an elementary $i$-characteristic $C_0'$ which joins two points of $B$ both of which are within $2\epsilon_0$ of the left endpoint of $B$. We work with the right half $B_R$ of $B$; the endpoints of $B_R$ are 0 and a point $d$ to the right of 0, with $|d| = 4$. In light of our assumption (4.2), it follows from the limit characteristic principle that there are two elementary $i$-arcs $F_1 = a_1 a_2$ and $F_2 = b_1 b_2$ joining points of $B_R$ such that
$$|a_1 - a_2|, |b_1 - b_2| < \tfrac{1}{H_0}$$
and
$$\tfrac{1}{H_0} < \text{dist}(F_1, 0), \text{dist}(F_2, d) < \tfrac{2}{H_0}.$$

Consider the characteristic subdomain $(D, C_0')$ for which bot$(D)$ contains the $\partial X \backslash B$, so that the only part of $X$ not in $D$ the union of $C_0'$ and the interior of the (small) subdomain bounded by $C_0'$ and the arc of $B$ which joins its endpoints. *For the rest of this step all order related statements refer to the order on* bot$(D)$. We apply Proposition 3.17 to $\min_D(F_1 \cup F_2)$ with $A = F_1$ and $B = F_2$ and examine in turn each of the three cases of its conclusion.

**Case (*i*).** If there were a contact point $p$ of $\min_D(F_1 \cup F_2)$ between $F_1$ and $F_2$ for which
$$\text{dist}(F_1, p), \text{dist}(F_2, p) > \tfrac{1}{H_0},$$
then Proposition 3.33 would imply that there is an ESF-$i$ subarc of $B_R$ of diameter $\epsilon_0$, contrary to (4.2). Thus, for every contact point $p$ of $\min_D(F_1 \cup F_2)$ between $F_1$ and $F_2$ one of
$$\text{dist}(F_1, p) \leq \tfrac{1}{H_0} \text{ or } \text{dist}(F_1, p) \leq \tfrac{1}{H_0}$$
must hold. Say for definiteness that the former holds for some contact point $p$ between $F_1$ and $F_2$, and let $\overline{a}$ be the rightmost such point. Then it is clear that $\min_D(F_1 \cup F_2)$ has an elementary subcharacteristic $C_0$ joining $\overline{a} \in B$ to a point $\overline{b} \in$ bot$(D)$, where $\overline{b}$ is either between $F_1$ and $F_2$ but within $\frac{1}{H_0}$ of $F_2$ or to the right of $F_2$ (with respect to the order on bot($D$)). In either event, $|\bar{b} - \bar{a}| \geq 1$. The case dist($F_2, p$) $\leq \frac{1}{H_0}$ is handled similarly.

For the other two cases we let $a$, $b$, $e$, $f$, $U$ and $L$ be as in the conclusion of Proposition 3.17. Furthermore, let the left and right endpoints of $C'_0$ be $l$ and $r$ (with respect to the order on bot($D$), so that since bot($D$) is very close to a segment of $\mathbb{R}$, $l$ lies to the right of $r$ in the usual sense).

Case (*ii*) Let $\alpha$ and $\beta$ be the left and right endpoints of $L$. If $|\beta - \alpha| \geq 1$, then we can take $L$ for $C_0$ since $\alpha$ is to the left of $F_1$ and therefore is in $B$. Thus we may assume that $|\beta - \alpha| < 1$. Starting at $a$ and traversing $\partial X$ in the *negative* direction the indicated points occur in the following order: $a$, $\beta$, $a_2$, $a_1$, $\alpha$, $e$, $l$, $r$, $f, b$ (with some equalities possible). If $|e - f| < 1$ then the arc $e\alpha$ joins $e$, $\alpha \in B$ and $|e - \alpha| \geq 2$. It is then easy to see that either $e\alpha$ has a proper contact point $g$ for which $|g - e|$, $|g - \alpha| > \frac{1}{H_0}$, which (in light of Proposition 3.33) would violate (4.2), or $e\alpha$ has a subarc $C_0$ which is an elementary $i$-characteristic with the desired properties. Thus we only need consider what happens if $|e - f| \geq 1$. But in this case $C_0 = ef$ has the desired properties since it is an elementary $i$-characteristic and $e \in B$.

Case (*iii*) Here we assume that we are not simultaneously in Case (*ii*), since were that so we would be done. The point $b$ is obviously between $F_1$ and $F_2$, but, since we are not in Case (*ii*), $a$ is between $l$ and $F_1$. Let $L = pq$, with $p < q$. If $|b - p| \geq \frac{1}{H_0}$, then again by Proposition 3.33, $B$ would have an ESF-$i$ arc of diameter $\epsilon_0$ contradicting (4.2). Thus we can assume that $|b - p| < \frac{1}{H_0}$. If $|p - q| < 1$, then we will have $|b - a| > 1$ and we can take for $C_0$ any elementary $i$-characteristic closely approximating the extended characteristic $ab$ (recall the construction of extended characteristics in Section 3), since its endpoints can be made arbitrarily close to $a$ and $b$, and consequently can be chosen to lie in $B$. On the other hand, if the distance between $p$ and $q$ is at least 1, then we can take $C_0$ to be the elementary $i$-characteristic $pq$, since $p \in B$. Thus, we are done with Step 1.

Step 1 says that there is a subarc $P$ of $\partial X$ and an elementary $i$-characteristic $C_0$ joining the two endpoints of $P$ such that $P$ contains a subarc $I$ of $B$ the distance between whose endpoints is at least 1 and one of whose endpoints is an endpoint of $C_0$. To be specific we assume that the left endpoint of $I$ is the one $I$ has in common with $C_0$. We shall henceforth work with the characteristic subdomain $(D_0, C_0)$ bounded by $P \cup C_0$, so that in particular all order related statements refer to the order on bot($D_0$) unless otherwise indicated. We give $P$ the usual arc length parametrization $z = z(s)$, $0 \leq s \leq L$, with increasing $s$ corresponding to the positive orientation on $\partial X$ and such that $z(0)$ is the left endpoint of $I$. Obviously, $z([0,1]) \subset I$. By the bound of $\frac{1}{10000}$ on the unsigned curvature of $B$ it is clear that

(4.3) $\qquad |z(s_2) - z(s_1)| \geq \frac{9}{10}|s_2 - s_1|$, for $z(s_1)$, $z(s_2) \in I$.

(Obviously, a constant considerably closer to 1 than $\frac{9}{10}$ would also work here.) We also observe that for $\alpha$, $\beta \in \text{bot}(D_0)$

(4.4) $\qquad |\beta - \alpha| \geq \gamma_0 \min\{1, \lambda([\alpha', \beta'])\}$, for all $[\alpha', \beta'] \subset [\alpha, \beta]$.

This follows from (4.1) and the fact that one of the arcs of $\partial X$ joining $\alpha$ to $\beta$ contains $[\alpha', \beta']$ and the other contains the endpoints of $C_0$, the distance between which is at least 1.

**Step 2**. There is a subarc $J = pq$ of $I$ with $|q - p| \geq \frac{1}{40}$ such that $p$ and $q$ are joined by an elementary $i$-characteristic in $D_0$.

To establish this we again we use an argument based on the case that results when Proposition 3.17 is applied. However, for this step this proposition must be applied, in conjunction with the trapped area bound (Propositions 3.34), in a sequential manner. By our underlying assumption (4.2) together with the limit characteristic principle there are elementary $i$-characteristics $A_1, \ldots, A_{H_0} \in \mathcal{I}(D_0)$ for which the arc $I_k \subset I \subset \text{bot}(D_0)$ joining the endpoints of $A_k$ lies in the middle third of $z([\frac{k-1}{H_0}, \frac{k}{H_0}])$ (i.e., lies in $z([\frac{k-\frac{2}{3}}{H_0}, \frac{k-\frac{1}{3}}{H_0}])$), $1 \leq k \leq H_0$. Note that by (4.3)
$$(4.5) \qquad \text{dist}(I_{k_1}, I_{k_2}) \geq \tfrac{9}{10H_0}(|k_2 - k_1| - \tfrac{1}{3})$$

Assume inductively that we have elementary $i$-characteristics $C_k$, $0 \leq k \leq t$, such that $C_k$ envelopes only $A_{l(k)}, A_{l(k)+1}, \ldots, A_{r(k)}$, but none of the other $A_n$, with
$$C_0 \succeq C_1 \succeq \ldots \succeq C_t,$$
$$1 = l(0) \leq l(1) \leq \ldots \leq l(t) < r(t) \leq r(t-1) \leq \ldots \leq r(0) = H_0$$
and for which $r(k) - l(k)$ is strictly decreasing and satisfies
$$(4.6) \qquad r(k) - l(k) \geq \tfrac{H_0}{3}, 0 \leq k \leq t.$$
We apply Proposition 3.17 with $A = A_{l(t)}$ and $B = A_{r(t)}$ to analyze $C = \min_{D_0}(A_{l(t)} \cup A_{r(t)})$ and consider separately each of the three cases in the conclusion of that proposition. More specifically, we will show that as long as Case (*iii*) does not occur we will have obtained the desired $J$ or produce $C_{t+1}$ with
$$(4.7) \qquad r(t+1) - l(t+1) = r(t) - l(t) - 1,$$
that Case (*iii*) will have to occur long before (4.6) can be violated, and that when Case (*iii*) occurs we have the desired $J$.

Case (*i*) There can be no contact points of $C$ between $A_{l(t)+1}$ and $A_{r(t)-1}$, since were there to be such a point $p$, the distance from $p$ to each of the endpoints of $C$ would, by (4.5), have to be at least $\frac{9}{10H_0}(\frac{2}{3}) > \frac{1}{2H_0}$, so that we can apply Proposition 3.33 to conclude that there is an ESF-$i$ arc of diameter $\epsilon_0$ on $B$, which contradicts (4.2). If there is a contact point of $C$ between $A_{r(t)-1}$ and $A_{r(t)}$ we will have obtained the desired conclusion because in that case $C$ would have to have a subarc $E$ which is an elementary $i$-characteristic with left endpoint to the left of $A_{l(t)+1}$ and right endpoint between $A_{r(t)-1}$ and $A_{r(t)}$, so we stop the process here, Step 2 having been established with with $J = E$ in light of (4.6). If there is no contact point of $C$ between $A_{r(t)-1}$ and $A_{r(t)}$, then there is a contact point of $C$ between $A_{l(t)}$ and $A_{l(t)+1}$ and $C$ has a subarc $C_{t+1}$ which is an elementary characteristic whose left endpoint is between $A_{l(t)}$ and $A_{l(t)+1}$ and whose right endpoint lies to the right of $A_{r(t)}$. We have $l(t+1) = l(t) + 1$ and $r(t+1) = r(t)$, so that (4.7) holds.

For the other two cases we use the notation of Proposition 3.17: $L$, $U = C = \min_{D_0}(A_{l(t)} \cup A_{r(t)}) = ef$, and $ab$ is the nonmonotone extended $i$-characteristic.

Case (*ii*) but not Case (*iii*). Because we are not in Case (*iii*), $b$ lies to the right of $A_{r(t)}$. We have that $a$ is between $A_{l(t)}$ and $A_{l(t)+1}$, since otherwise $|e - a| \geq \frac{9}{10H_0}(\frac{2}{3}) > \frac{1}{2H_0}$, by (4.5), and then, since $[e, a] \subset [e, b]$,
$$|e - b| \geq \gamma_0 \min\{1, \lambda([e, a])\} \geq \gamma_0 \min\{1, |e - a|\} > \tfrac{\gamma_0}{2H_0}$$

by (4.4), so that by Proposition 3.33 applied to the contact point $e$, there is an ESF-$i$ of diameter $\epsilon_0$ on $B$ which contradicts (4.2). From the definition of extended characteristics there are elementary $i$-characteristics $U' = a'b' \preceq U$ with $a'$ and $b'$ arbitrarily close to $a$ and $b$, respectively. The point $a'$ can therefore be taken to lie between $A_{l(t)}$ and $A_{l(t)+1}$ and $b'$ can be taken to lie to the right of $A_{r(t)}$. In this case we let $C_{t+1}$ be any such $U'$ and set $l(t+1) = l(t)+1$ and $r(t+1) = r(t)$, so that we again have (4.7).

Case (*iii*). We will show that the first time this case occurs desired conclusion of Step 2 holds. Thus we assume for all $t < t_0$ we were in either Case (*i*) or Case (*ii*). Now we can only say that $a$ lies somewhere to the right of the left endpoint of $A_{l(t_0)}$ and that $f$ lies somewhere to the left of $A_{r(t_0)}$, but since we cannot use the Proposition 3.33 to any effect (because $f$ is not necessarily in $B$), we cannot conclude, in analogy with the preceding case, that $b$ lies between $A_{r(t_0)-1}$ and $A_{r(t_0)}$. If $f \in B$, then by (4.3) and (4.6) it follows that $|f - e| \geq \frac{1}{6}$. But $f$ is a point of $\partial X$, and if it is not in $B$, then (refer to the definition of the bottom $B$ of $X$) it is certainly not $N(e,1) \cap B$, so that by property (4b) of $X$, we have $|f - e| \geq \frac{1}{2}$. Therefore, no matter what, $|f - e| \geq \frac{1}{6}$. Let $D'$ be the subdomain of $D_0$ bounded by $ab \cup [a,b]$. It follows immediately from (4.3) that $\lambda([a,b]) \leq \frac{10}{9}|b-a|$. From the definition of extended characteristics it follows that there elementary $i$-characteristics $E$ inside $D'$, joining points of $a', b' \in [a,b]$ which are arbitrarily close to $a$ and $b$, respectively. But it follows as a simple corollary of Proposition 3.6 that
$$\tfrac{1}{6} \leq |f - e| \leq \text{diam}\,(ef) \leq \text{diam}(ab) \leq 5\lambda([a,b]) \leq \tfrac{50}{9}|b-a|,$$
so that, $|b - a| > \frac{1}{40}$. But taking $|a' - a|$ and $|b' - b|$ sufficiently small, the corresponding $E$ will serve as the desired $J$.

Thus either we stop at $t$ with the desired $J$, or go on to $t + 1$, the latter occurring only when we are in Case (*i*) or Case (*ii*) but not Case (*iii*). We now show we must actually arrive at Case (*iii*) long before (4.6) can be violated. Say that we have arrived at $t = T < \frac{H_0}{3}$. Let $T' = [\frac{T}{2}]$. At least one of the following two things must have occurred:
(A) At least $T'$ of the times that we passed from $t$ to $t+1$ we will have done so because we are in Case (*i*) there is a proper contact point $p_t$ of $\min_{D_0}(A_{l(t)} \cup A_{r(t)})$ between $A_{l(t)}$ and $A_{l(t)+1}$ and $p_t$ is the left endpoint of $C_{t+1}$, or
(B) At least $T'$ of the times that we passed from $t$ to $t+1$ we will have done so because we are in Case (*ii*) but not Case (*iii*).

We deal with possibility (A) first. Say $k_1 < k_2$ are two values of $t$ for which we are in Case (*i*). Let $C'$ be the elementary $j$-characteristic one of whose endpoints is $p_{k_2}$. Let its second endpoint be $p'$. If $C'$ were to cross $C_{k_1+1}$, then, since a $j$-characteristic can have at most one point in common with an $i$-characteristic, the other endpoint of $C'$ would lie to the left of $A_{l(k_1)}$ or to the right of $A_{r(k_1)}$, and the distance between its endpoints would, by (4.5), have to be at least $\frac{9}{10H_0}(\frac{2}{3}) > \frac{1}{2H_0}$. However, (as we have seen happen before) by (4.6) together with Proposition 3.33, $B$ would then have an ESF-$i$ arc of diameter $\epsilon_0$, contradicting (4.2). Thus $C'$ does not cross $C_{k_1+1}$, so that $p' \in \text{bot}(D_0)$. If $\alpha$, $\beta$ are points of $\text{bot}(D_0)$ with $\alpha$ lying to the left of $C_{k_2+1}$ and $\beta$ to its right, then $\lambda(\alpha\beta) \geq \frac{1}{4}$ since $T < \frac{H_0}{3}$, so that by (4.4),
(4.8) $$|\beta - \alpha| \geq \tfrac{\gamma_0}{4}.$$

If $p'$ lies to the right of $C_{k_2+1}$, then each point of $C_{k_2+1} \cap X$ is joined to a point in bot$(D_0)$ to the right of $C_{k_2+1}$ by a $j$-half-characteristic lying between $C_{k_2+1}$ and $C' \preceq C_{k_1+1}$. Then a simple argument (which was used in the proof of Proposition 3.34) based on Proposition 2.8 and the fact that for each point of the first half of $C_{k_2+1}$, which as length at least $\frac{\gamma_0}{8}$, the corresponding $j$-half-characteristic has length at least $\frac{\gamma_0}{8}$, shows that the area between $C_{k_2+1}$ and $C_{k_1+1}$ must be at least $\eta(\frac{\gamma_0}{4})^2 = \frac{\eta\gamma_0^2}{16} \geq \frac{\eta_1\gamma_0^2}{64}$. If, on the other hand, $p' < p_{k_2}$, then $|\Delta\theta(C')| \geq \frac{\pi}{4}$. Also, as we just saw (4.8) holds for all $\alpha$, $\beta \in$ bot$(D_0)$ with $\alpha$ lying to the left of $C_{k_2+1}$ and $\beta$ to its right. Thus by the trapped area bound (Proposition 3.34) it follows that the area between $C_{k_2+1}$ and $C_{k_1+1}$ is at least $\frac{\eta_1\gamma_0^2}{16} > \frac{\eta_1\gamma_0^2}{64}$, so that this lower bound holds no matter which side of $C_{k_2+1}$ the point $p'$ lies on. But then $\frac{(T'-1)\eta_1\gamma_0^2}{64} \leq \mu(X) \leq 200$, so that $T' \leq \frac{12800}{\gamma_0^2 \eta_1} + 1$. This in turn means that $T \leq \frac{30000}{\gamma_0^2 \eta_1} < \frac{H_0}{10}$ (since $\gamma_0$ and $\eta_1$ are both in $(0,1)$).

The other possibility (B) is handled in a similar manner. Say $k_1 < k_2$ are two values of $t$ for which we are in Case (*ii*) but not in Case (*iii*) when passing from $C_t$ to $C_{t+1}$. From the construction of $C_{k_1+1}$ we have that $U = \min_{D_0}(A_{l(k_1)} \cup A_{r(k_1)}) = u_1 u_2$ is a subarc of a nonmonotone extended $i$-characteristic $W$ which has another subarc $S$ which is also an elementary $i$-characteristic for which $S \preceq U$ and $A_{l(k_1)} \preceq S$. In addition $S$ precedes $C_{k_1+1} \preceq U$, so that in particular $S$ precedes $A_k$ for all $k > l(k_1)$. There is an elementary $j$-characteristic $S'$ which joins $u_1$ to a point $u' \in [u_1, u_2]$. It follows immediately from (4.6) and (4.4) that $|u_2 - u_1| \geq \frac{\gamma_0}{4}$. We claim that $S' \cap D_0 \cap C_{k_2+1} = \emptyset$. To see this, assume to the contrary that $z \in S' \cap D_0 \cap C_{k_2+1}$. It then follows from (4.4) and (4.5) that $|u' - u_1| \geq \frac{9\gamma_0}{10H_0}(\frac{2}{3}) > \frac{\gamma_0}{2H_0}$, so that by Proposition 3.33, $B$ would have an ESF-$i$ arc of diameter a least $\epsilon_0$ in contradiction to our underlying assumption (4.2). Now we apply the trapped area bound to the $i$-characteristic subdomain $(D_1, C_{k_2+1})$, where $D_1$ is bounded by the curve made up of $C_{k_2+1}$ together with $\partial X \setminus (\alpha_1, \alpha_2)$, where $\alpha_1 < \alpha_2$ are the endpoints of $C_{k_2+1}$ (that is, $D_1$ is the part of $X$ that remains when the part of $D_0$ on and below $C_{k_2+1}$ is removed). In particular the interior of $S'$ lies in $D_1$ and $U \in \mathcal{I}(D_1)$. It is also clear from (4.6) and (4.4) that (4.8) holds for $\alpha$ and $\beta$ in bot$(D_1)$ on opposite sides of $C_{k_2+1}$. For the same reason that $S' \cap D_0 \cap C_{k_2+1} = \emptyset$, we have that $u'$ must lie between $u$ and $\alpha_1$ with respect to the order on bot$(D_0)$. This means that both endpoints of $S'$ are in $B$ so that we clearly have $\Delta\theta(S') \geq \frac{\pi}{4}$. It also means that with respect to the order on bot$(D_1)$, $u_1$ lies between $u'$ and $u_2$. We can now apply the trapped area proposition (with $(D_1, C_{k_2+1})$ playing the role of $(D, V)$) to deduce that the area between $U$ and $C_{k_2+1}$ is bounded below by $\frac{\eta_1\gamma_0^2}{16}$, so that this same lower bound holds for the area between $C_{k_1}$ and $C_{k_2+1}$. If the values if $t$ in question are $t_1 < t_2 < \ldots t_{T'}$, then if $R_k$ is the region between $C_{t_k}$ and $C_{t_{k+1}+1}$, then $R_1, R_3, \ldots$ are disjoint so that

$$\left(\frac{T'-1}{2}\right)\frac{\eta_1\gamma_0^2}{16} \leq \mu(X) \leq 200,$$

and therefore $T' \leq \frac{6400}{\gamma_0^2 \eta_1} + 1$. Thus as in the case that (A) holds we again have $T < \frac{H_0}{10}$. This concludes the proof of Step 2.

**Step 3**. Now we work with the $J = pq \subset I$, with $|q - p| \geq \frac{1}{40}$, whose endpoints are joined by a $C = pq \in \mathcal{I}(D_0)$. Clearly, there are two points $a_1, a_2 \in J$ such that if we consider $X_k = \Sigma(a_k, \frac{1}{10000})$, $k = 1, 2$, where the $\Sigma$ (as defined in the second paragraph of this section) are with respect to the normalized domain $X$ with which we are now working, then

$10000(X_k - a_k)$ is a normalized domain, $k = 1, 2$,
dist$(X_k, \{p, q\}) \geq \frac{1}{500}$, $k = 1, 2$,
dist$(X_1, X_2) \geq \frac{1}{500}$, and
the bottoms of $X_1$ and $X_2$ are in $J$.

By what we have shown it then follows that either the bottom one of $X_1$ or $X_2$ has an ESF-$i$ of length $\frac{\epsilon_0}{10000}$, in which case we have reached the desired conclusion, or for both $k = 1$ and $k = 2$ the bottom of $X_k$ has a subarc $p_k q_k$, with

(4.9) $$|q_k - p_k| \geq \frac{1}{400000}$$

whose endpoints are joined by an elementary $i$-characteristic $C_k$. In this latter case, since $C$ joins the endpoints of $J$, $C_k \preceq C$, $k = 1, 2$, we can apply Proposition 3.17 to $\max_{D'}(C_1 \cup C_2)$, where $D'$ is the characteristic subdomain bounded by $C \cup [p, q]$. But in light of the above conditions satisfied by the $X_k$ and (4.9), and since all contact points of $\max_{D'}(C_1 \cup C_2)$ are in $J$, Proposition 3.33 immediately implies that there is a $\gamma_1 = \gamma_1(K)$ for which $J$ has an ESF-$i$ arc of $J$ of diameter $\gamma_1$. Thus we have proved that the bottom of a normalized domain has an ESF-$i$ arc of length $\Delta_0 = \min\{\frac{\epsilon_0}{10000}, \gamma_1\}$.∎

**Proof of Main Theorem**. It follows from the opening discussion and the proposition we have just proved, that there are $\delta_0 = \delta_0(G)$, $\alpha_0 = \alpha_0(K) \in (0, 1)$ such that any arc $B$ of $\partial G$ of diameter at most $\delta_0$ has an ESF-$i$ subarc $B'$, for which diam$(B') > \alpha_0$diam$(B)$. It is then clear that there numbers $\xi_0 = \xi_0(K)$, $\tau_0 = \tau_0(K)$, both in $(0, 1)$, such that every arc $B$ with diam$(B) \leq \delta_0$ has two disjoint subarcs $B_1$ and $B_2$ for which

diam$(B_1)^{\tau_0}$ + diam$(B_2)^{\tau_0}$ < $\xi_0$diam$(B)^{\tau_0}$,
diam$(B_1)$, diam$(B_2)$ < $(1 - \alpha_0)$diam$(B)$
$B\setminus(B_1 \cup B_2)$ is ESF-$i$.

Let diam$(B) \leq \delta_0$. We start with two such arcs $B_1$ and $B_2$ corresponding to $B$, then we apply the same fact to each of these to get four arcs of diameter $l_1, \cdots l_4 < (1 - \alpha_0)^2$ for which $\sum l_k^{\tau_0} < \xi_0^2 \delta_0^{\tau_0}$, and such that the complement of their union is ESF-$i$, then we do so again to get eight arcs of length at most $(1 - \alpha_0)^3$ for which the corresponding $\sum l_k^{\tau_0}$ is less than $\xi_0^3 \delta_0^{\tau_0}$ and for which the complement of their union is ESF-$i$, and so on. At the $n^{\text{th}}$ stage we have $2^n$ arcs $B_k^{(n)}$ for which $\sum_k (\text{diam}(B_k^{(n)}))^{\tau_0} < \xi_0^n \delta_0^{\tau_0}$, which tends to 0, and for which $B\setminus(\underset{k}{\cup} B_k^{(n)})$ is ESF-$i$. At this stage the diameter of the largest of the $2^n$ arcs is at most $(1 - \alpha_0)^n$. But then $B$ is the union of a set of $\tau_0$-dimensional Hausdorff measure 0 and a set that has at most countably many $i$-singularities. In light of Corollary 3.29 and the fact that there are at most countably many vertices of fans on $\partial G$ it follows that the set of singularities on $B$, and therefore on all of $G$ has Hausdorff dimension at most $\tau_0$.∎

**Corollary 4.2.** *Let $\Theta$ be a normal system. Then there is a $\tau = \tau(\Theta) < 1$ such that for any smoothly bounded $G \subset \mathbb{C}$ and any locally Lipschitz solution $R$ of $\Theta$ on $G$ the set of points of $\partial G$ at which $R$ does not have a nontangential limit has Hausdorff dimension at most $\tau$.*

**Proof**. This is an immediate consequence of the Main Theorem and Propositions 3.21 and 3.26.∎

## 5. Construction of solutions with singularity sets of positive Hausdorff dimension

Although we have chosen to carry out our construction in $\mathbb{H}$ to avoid cumbersome arguments, suitable, largely straightforward changes will allow an analogous procedure to be carried out in any smoothly bounded domain. Throughout this section $\Theta$ will denote any fixed normal system. We begin with a brief discussion of characteristic initial value problems for $\Theta$, since our construction largely proceeds by joining together solutions of appropriate instances of such problems in domains sharing a boundary characteristic.

Let $C_k$ be $C^\infty$ curves parametrized by $z_k(s)$, $0 \leq s \leq l_k$, $k = 1, 2$ for which $z_i(0) = z_j(0)$ and for which $z_i'(0)$ and $z_j'(0)$ are mutually orthogonal. (In fact $l_k$ can be infinite, as is the case, for example, when $C_k$ is a ray.) We want to construct a solution of the normal system $\Theta$ for which $C_k$ is an $k$-characteristic arc, $k = 1, 2$. From the definition of normal system, for any $\rho_1$ there is a discrete set of values $\rho_2$ such that $e^{i\theta(\rho_1,\rho_2)}$ and $ie^{i\theta(\rho_1,\rho_2)}$ are tangent at $z_1(0)$ to $C_1$ and $C_2$, respectively. Given any such *admissible corner value* $\rho = (\rho_1, \rho_2)$ there are unique continuous functions $R_k(s)$ such that $e^{i\theta(\rho_1,R_2(s))}$ and $ie^{i\theta(R_1(s),\rho_2)}$ are tangent to $C_1$ and $C_2$ at the points $z_1(s)$ and $z_2(s)$, respectively. For $t = (t_1, t_2)$ we define
$$\overline{R}(t) = (R_1(t_2), R_2(t_1)), 0 \leq t_k \leq l_k, k = 1, 2$$
and
$$\overline{\theta}(t) = \theta(R_1(t_2), R_2(t_1)), 0 \leq t_k \leq l_k, \ k = 1, 2.$$
We seek $\zeta : [0, l_1] \times [0, l_2] \to \mathbb{C}$, where $\zeta(t) = \xi(t) + i\eta(t)$, for which $e^{i\overline{\theta}(t)}$ and $ie^{i\overline{\theta}(t)}$ are tangent to the curves $\zeta([0, l_1], t_2)$ and $\zeta(t_1, [0, l_2])$, respectively at the point $\zeta(t)$ and which satisfies the initial conditions
(5.1) $$\zeta(t_1, 0) = z_1(t_1) \quad \text{and} \quad \zeta(t_2, 0) = z_2(t_2).$$
The tangency condition can be written as the linear hyperbolic system
(5.2) $$\cos\overline{\theta}\,\tfrac{\partial\eta}{\partial t_1} - \sin\overline{\theta}\,\tfrac{\partial\xi}{\partial t_1} = 0; \quad \sin\overline{\theta}\,\tfrac{\partial\eta}{\partial t_2} + \cos\overline{\theta}\,\tfrac{\partial\xi}{\partial t_2} = 0,$$
The problem (5.2) with initial conditions (5.1) is well posed and has a $C^\infty$ solution $\overline{\theta}$ for any $C^\infty$ initial curves $C_1$ and $C_2$. First consider the case in which the $\zeta(t) = \zeta(i, C_i, C_j, \rho, t)$ determined in this manner has an everywhere nonzero Jacobian determinant and is globally one-to-one. We then have a solution to the system $\Theta$ in $\zeta([0, l_1] \times [0, l_2])$, which is a characteristic quadrilateral, given by
(5.3) $$R(\zeta(t)) = \overline{R}(t).$$
We refer to this solution as $R(i, C_i, C_j, \rho, z)$ and to the characteristic quadrilateral in which it is defined as $Q(i, C_i, C_j, \rho)$. This is just a mechanism for indicating which of the two curves is the $i$-arc. Obviously,
$$R(j, C_j, C_i, \rho, z) = R(i, C_i, C_j, \rho, z) \text{ and } Q(j, C_j, C_i, \rho) = Q(i, C_i, C_j, \rho),$$
and analogously for $\zeta$. Even if $\zeta$ is not one-to-one on $[0, l_1] \times [0, l_2]$, then (5.3) gives a multivalued solution of $\Theta$ on any domain $\zeta(E)$, provided that $\zeta$ is a local diffeomorphism on $E$. Let $i$ and $j$ be such that $\arg\{z'_j(0)\} = \arg\{z'_i(0)\} + \frac{\pi}{2}$. It is a well-known property of the characteristic initial value problem for genuinely nonlinear systems that if $\frac{d\arg\{z'_i(t)\}}{dt} \leq 0$ on $[0, l_1]$ and $\frac{d\arg\{z'_j(s)\}}{ds} \geq 0$ on $[0, l_2]$, then $\zeta$ will be locally diffeomorphic on $[0, l_1] \times [0, l_2]$. This is simply a reflection of the fact that, in light of the quasi-HP nature of the the associated inclination function and the length monotonicity property (Proposition 2.7), under these hypotheses the curvature of $j$-characteristics will not blow up as one moves along $i$-characteristics away from the convex side of a $j$-characteristic $C$, and analogously when the roles of $i$ and $j$ are interchanged. We note that this property holds in the particular case in which one of the initial characteristics is a line segment or ray. We also note that it holds even if one or both of the initial curves is not a simple arc, so that $\zeta$ will be a local homeomorphism if, for example, one of them is a circle covered several times.

We next need to discuss briefly the smooth adjunction of line segments and circles to $C^\infty$ curves in a specific, constructive fashion; this is another essential element of the construction process. Let $C$ be a $C^\infty$ curve parqametrized by $Z(s)$, $0 \leq s \leq l + \delta$. We want to perform the adjunction with no change in $Z$ on $[0, l]$. Let $\kappa_0$ be any number. It is clear that there is an operator $F(Z,l,\delta,\kappa_0,\tau)$ such that $w = F(Z,l,\delta,\kappa_0,\tau) \in C^\infty([0, \infty))$ has the following properties:
(i) $w(s) = Z(s), 0 \leq s \leq l$,
(ii) If $w'(s) = e^{i\phi(s)}$, then $\phi(s) = (s - l - \delta)\kappa_0 + \tau$, $s \geq l + \delta$.
It is clear how this can be done. Indeed, if $Z'(s) = e^{i\alpha(s)}$ on $[0, l + \delta]$ we let $\beta(s)$ be the continuous function which coincides with $\alpha(s)$ on $[0, l + \delta/3]$, which is given by $(s - l - \delta)\kappa_0 + \tau$ for $s \geq l + 2\delta/3$ and is linear on $[l + \delta/3, l + 2\delta/3]$. We take for $\sigma(s)$ a specific nonnegative $C^\infty$ function on $\mathbb{R}$ with support in $(-1, 1)$ and $\int_{-\infty}^{\infty} \sigma(s)ds = 1$ and then convolve $\beta(s)$ with $\frac{10}{\delta}\sigma(\frac{10s}{\delta})$. The desired $w = F(Z,l,\delta,\kappa_0,\tau)$ is defined by with $w'(s) = e^{i\beta(s)}$, $w(0) = Z(0)$. We note that if $C$ is convex to the right (left) and $\kappa_0$, $\tau - \kappa_0\delta/3 - \alpha(l + \delta/3) \geq 0$ ($\leq 0$), then the curve given by $F(Z,l,\delta,\kappa_0,\tau)$ is concave towards the same side as $C$.

We now introduce notation and terminology to be used in our construction. We denote by $\mathcal{CL}$ the family of $C^\infty$ arcs with initial and terminal straight subarcs. Let $\mathcal{S}(\epsilon)$ denote the class of $C^\infty$ curves parametrized by $z = z(s)$, $0 \leq s \leq L$ with the following properties.
(i) $\Im\{z(0)\} = \Im\{z(L)\} = -\epsilon$,
(ii) $\arg\{z'(0)\} = \frac{\pi}{2} + \epsilon$ and $\arg\{z'(L)\} = -(\frac{\pi}{2} + \epsilon)$,
(iii) $\frac{d\arg\{z'(s)\}}{ds} \leq 0$, $0 \leq s \leq L$
(iv) $C \subset N(N(0,1) \cap \mathbb{H}, 2\epsilon)$
(v) $C \in \mathcal{CL}$
For $t > 0$ and $\alpha \in \mathbb{R}$ we denote $\{tS + \alpha : S \in \mathcal{S}(\epsilon)\}$ by $\mathcal{S}(\epsilon, t, \alpha)$. If $C$ is any arc joining $a < b$ in $\overline{\mathbb{H}}$, $B(C)$ will denote the set whose boundary is $C \cup [a, b]$. For two such arcs $C_1$, $C_2$ with $C_1 \subset B(C_2)$ we denote the closure of $B(C_2) \backslash B(C_1)$ by $E(C_1, C_2)$. In what follows we will often deal with a $C^\infty$ solution $R$ defined only in a neighborhood in $B(C)$ of $C$, where $C$ is an $i$-arc of $R$ (by which we mean that $R$ can be extended to a $C^\infty$ solution in a two-sided neighborhood of $C$ and that $C$ is an $i$-arc of this extension). In this case we shall call $C$ a *one-sided characteristic* of $R$. An $i$-characteristic $C$ of a solution $R$ will be said to satisfy the *orthogonal segment condition* (OSC) if the $j$-characteristic passing through each point $c \in C$ all contain a straight line subarc containing $c$ in its interior; obviously it is enough for a single point $c \in C$ to have this property for the OSC to hold. We extend the use of the term OSC in the obvious manner to apply to one-sided characteristics. If $D$ is a domain separated into two subdomains $D_1$ and $D_2$ by an arc $C$, and $R^{(1)}$ and $R^{(2)}$ are $C^\infty$ solutions in these domains, respectively, for which $C$ is a one-sided $i$-characteristic satisfying the OSC, then together $R^{(1)}$ and $R^{(2)}$ give a $C^\infty$ in all of $D$. Our construction uses this simple fact, which allows for the smooth pasting together of solutions of characteristic initial value problems in contiguous domains.

Next we discuss a specific class of characteristic initial value problems in which one of the initial characteristics is a ray. Let $C$ be parametrized by $z(s)$ satisfying

(5.4) $\qquad \frac{d\arg\{z'(s)\}}{ds} \leq 0, 0 \leq s \leq L$ with $\Im\{z(0)\}, \Im\{z(L)\} < 0$,

and

(5.5) $\qquad \arg\{z'(0)\} \geq \frac{\pi}{2} + \epsilon$ and $\arg\{z'(L)\} \leq -(\frac{\pi}{2} + \epsilon)$

and consider the solution of the initial characteristic value problem where $C_i = C$ and $C_j$ is a ray orthogonal to $C$ at $z(0)$ and emanating to the left of $C$. Then the $j$-characteristics of any corresponding solution are all rays orthogonal to $C_i$ and the $i$-characteristics are the orthogonal trajectories of this family of rays, and in fact are the curves $C(r)$, $r > 0$, parametrized by $z(t) + rz'(t)i$, $0 \leq t \leq L$. As $r$ tends to $\infty$ the $C(r)$ tend to arcs of a circle of radius $r$ and radian measure $\arg\{z'(0)\} - \arg\{z'(L)\}$.

We next construct two special $C^\infty$ arcs $U_i$. The arc $U_1$, parametrized by $u_1(s)$, $0 \leq s \leq \lambda_1$ will have the following properties.

(i) $\frac{d\arg\{u_1'(s)\}}{ds} \leq 0, 0 \leq s \leq \lambda_1$,

(ii) $\Im\{u_1(0)\} = -\epsilon$ and $u_1(\lambda_1) \in \mathbb{R}$

(iii) $U_1 \in \mathcal{CL}$,

(iv) $u_1'(0)\} = e^{i(\frac{\pi}{2}+\epsilon)}$,

(v) $u_1'(s) = e^{-i\pi/4}$, $\lambda_1 - d \leq s \leq \lambda_1$, for some $d > 0$,

(vi) $U_1 \cap N(0,2) = \emptyset$,

(vii) Let any $C'$ be any curve in $\mathcal{S}(\epsilon)$ and let $\alpha = \alpha(C)$ be such that the initial point of $C = C' + \alpha \in \mathcal{S}(\epsilon, 1, \alpha)$ is $p = -1 - \epsilon i$. Let $\rho \in \mathbb{R}^2$ be such that $ie^{i\theta(\rho)}$ is tangent to $C$ at its left endpoint $p$. Then there is a solution to the system $\Theta$ in $E(C) = E(C \cap \mathbb{H}, U_1 \cap \mathbb{H})$ for which $R(p) = \rho$, for which $C$ is a 2-characteristic and $U_1$ is an 1-characteristic, and for which both $C$ and $U_1$ satisfy the OSC.

We stress that the idea is that for the appropriate translate $C$ of *any* $C' \in \mathcal{S}(\epsilon)$ and any admissible value $\rho$ of $R$ at the initial point (i.e., left endpoint) $c_l$ of $C$ we have a solution in the domain $E(C \cap \mathbb{H}, U_1 \cap \mathbb{H})$. If we can construct such a curve it is obvious that it will have a right-hand counterpart $U_2$, having properties corresponding to (*i*)-(*vii*). In regard to (*v*), $U_2$ will terminate in a line segment of slope $\frac{\pi}{4}$ at its left end. In regard to (*vii*), $U_2$ and $C$ will be 2- and 1-characteristics, respectively, and the initial value $\rho$ will be the value of $R$ the left end of $C$, which in any case determines its value at the right endpoint of $C$ since the total curvature of $C$ is exactly $\pi + 2\epsilon$. For lack of a better term we shall refer to $U_1$ and $U_2$ as *universal* 1- *and* 2-*arcs* for the system $\Theta$. It is clear that translates $U_i + \alpha$, $\alpha \in \mathbb{R}$ of $U_i$ have an analogous universal property.

We show that such a $U_1$ exists, the existence of $U_2$ being identical apart from trivial details. Let $W$ be any $C^\infty$ arc parametrized by $w(s)$, $0 \leq s \leq L$ for which $\frac{d\arg\{w'(s)\}}{ds} \leq 0$ and which has total curvature $-\Delta < 0$. Let $V_1 \in \mathcal{CL}$ be parametrized by $v_1(s)$, $0 \leq s \leq 1$, with $v(0) = w(0)$, $v_1'(0) = iw'(0)$, for which $0 \geq \frac{d\arg\{v_1'(s)\}}{ds} \geq -\delta_0$, and whose total curvature is any $\delta_1 \leq \frac{\delta_0}{2}$. Let $A \leq |\frac{\partial \theta}{\partial R_k}| \leq B$ on $\mathbb{R}^2$ (as in the definition of normal system), so that the net of characteristics of any solution of the system $\Theta$ in a domain $G$ is a $K$-quasi-HP net on $G$ with $K = \frac{B}{A}$. Using this one can easily show that there is a $\delta_0 \in (0, \frac{\pi}{4})$ such that the corresponding characteristic coordinate mapping $\zeta(2, W, V_1, \rho, t)$ is a local diffeormophism on all of the corresponding characteristic coordinate rectangle $S = [0, L] \times [0, 1]$ for all such $W$ and $V_1$. To do so, one simply subdivides $S$ into small subrectangles and solves the corresponding initial value problem piece by piece. In light of the concavity of $V_1$ and the length change estimate (Proposition 2.4) the length of all of the 1-arcs parallel to $V_1$ in $\zeta(S)$ is at most 1. By the $K$-quasi-HP property all of the 2-arcs parallel to the initial 2-arc $W$ have total curvature at most $K\Delta$ so that from the length change estimate it follows that the length of the outer 2-characteristic, $\zeta(2, W, V_1, \rho, [0, L] \times \{1\})$, is at most $L + K\Delta$. Thus by a scale change we can smoothly add to $V_1$ another arc in $\mathcal{CL}$ of length $(L + K\Delta)/L$ whose parametrization $v_2(s)$ satisfies
$$0 \geq \tfrac{d\arg\{v_2'(s)\}}{ds} \geq -\delta_0 L/(L + K\Delta),$$
and whose total curvature is any $\delta_2 \leq \frac{\delta_0}{2}$ to obtain a longer 2-arc $V_2 \in \mathcal{CL}$, and accordingly extend the characteristic coordinate mapping to obtain $\zeta(2, W, V_2, \rho, t)$ which will still be locally one-to-one and will be $C^\infty$ in all of $[0, L] \times [0, L_2]$, where $L_2 = 1 + (L + K\Delta)/L$. Again by the length change estimate the outermost 2-arc will have length at most $L + L_2 K\Delta$. Clearly, we can continue this process so that we have $V_k$ made up of $k$ such arcs in $\mathcal{CL}$ of any total curvatures $\delta_1, \ldots, \delta_k \leq \frac{\delta_0}{2}$, where
$$\lambda(V_k) = \lambda(V_{k-1}) + (L + \lambda(V_{k-1})K\Delta)/L.$$
Applying this construction with $L$ and $\Delta$ chosen so as to accommodate all $W = C$, where $C$ is the appropriate translate of any element of $\mathcal{S}(\epsilon)$ as indicated in (*vii*), we see that there is some $M$ such that for some sufficiently small $\eta \in (0, \frac{\pi}{4})$ there is a curve $V$ parametrized by $v(s)$, $0 \leq s \leq M$ such that $V$ terminates in a straight line segment whose initial point is in $\mathbb{H}$ and whose slope is $-\tan \eta$. In fact we can allow the length of this segment to be as large as we want, so that we can take the imaginary part of the terminal point of $V$ to be $-\epsilon$. We next extend $V$ at its right end by smoothly adjoining a ray with slope exactly $-\frac{\pi}{4}$; (by means of $F(v, \lambda(V) - \frac{\epsilon}{2}, \frac{\epsilon}{2}, 0, -\frac{\pi}{4})$, as described above, where $v(s)$ parametrizes $V$); this extended curve will be called $V^1$. Note that for $\epsilon$ sufficiently small, the first subarc $V_1$ of $V^1$ will intersect $\mathbb{R}$ at a point $v_1(s)$ for which $\arg\{v_1'(s)\} \in (-\pi, -\pi + \frac{\delta_0}{2}) \subset (-\pi, -\frac{3\pi}{4})$, and will have the same convexity as $V$. We now consider the characteristic initial value problem with $V^1$ as the 1-arc and, as the 2-arc, a segment $A$ of length $r$ orthogonal to it at its initial point and emanating to the left of $V^1$. Because $V^1$ satisfies the OSC, the solution so generated will be $C^\infty$ in the domain made up of $E(C \cap \mathbb{H}, V^1 \cap \mathbb{H})$ together with translates of $A$ around $V^1$. Note that the corresponding solution is constant in any half-strip made up of all the translates of $A$ emanating from points on the ray in which $V^1$ ends. From this it is clear that for $r = \lambda(A)$ sufficiently large, an appropriate subarc of the translate $V^1(r)$ of $V^1$ to the outer endpoint of $A$ has properties (*i*), (*v*), (*vi*) and (*vii*). Properties (*ii*), (*iii*) and (*iv*), can be achieved in the following manner. First we take $r$ so large that at the left intersection point of $V^1(r)$ with $\mathbb{R}$ the angle formed by $V^1(r)$ is within $\frac{\epsilon}{4}$ of $\frac{\pi}{2}$. Then we remove the arc of $V^1(r)$ joining its initial point to a point with imaginary part between $-\frac{\epsilon}{2}$ and $0$ and at which the direction of the tangent vector is within $\frac{\epsilon}{2}$ of $\frac{\pi}{2}$. Finally, we can use the operator $F$ to smoothly adjoin an arc to the remaining part of $V^1(r)$ in such a way that (*ii*), (*iii*) and (*iv*) hold. The resulting curve is our $U_1$. Note that in (*vii*) we do not require the solution to bear any particular relation to the part of $U_1$ that lies in the lower half-plane.

We call the left and right endpoints of $U_1$, $e_1$ and $m_1$; the left and right endpoints of $U_2$ will be called $m_2$ and $e_2$ (*m* for middle and *e* for end, for reasons to be made clear momentarily). It is clear that $U_1$ depends solely on $\Theta$ and the parameter $\epsilon$ and that once these have been fixed, the solution $R$ generated in $E(C)$ depends on $C$ but that its values on $U_1$ depend only on $\rho = R(c_l)$, $c_l$ being the left endpoint of $C$, and analogously for $U_2$. For $\rho \in \mathbb{R}^2$ we denote by $f_k(\rho)$ the value of this solution $R$ at $m_k$. We also observe that in the case of $U_1$ we always have

(5.6) $\qquad \theta(f_1(\rho)) - \theta(\rho) = \theta(m_1) - \theta(c_l) = -(\frac{5\pi}{4} + \epsilon).$

As indicated above, in the case of $U_2$ we still take as $\rho$ the value of $R$ at the left endpoint $c_l$ of the initial curve $C$ used in the above construction which, as we have pointed out, uniquely determines the value of $R$ at the right endpoint of $C$ since the total curvature of all arcs in the family $\mathcal{S}(\epsilon)$ is $\pi + 2\epsilon$. Here we have that

(5.7) $\qquad \theta(f_2(\rho)) - \theta(\rho) = \theta(m_2) - \theta(c_l) = -(\pi + 2\epsilon) + (\frac{5\pi}{4} + \epsilon) = \frac{\pi}{4} - \epsilon.$

The $f_k$ are continuous functions of $\rho$ and there is a number $B_0$ (that only depends on the bounds on the $\frac{\partial \theta}{\partial R_k}$ associated with $\Theta$) such that

(5.8) $\qquad \|f_k(\rho) - \rho\| \leq B_0, \ \rho \in \mathbb{R}^2, k = 1, 2,$

where the norm is just the Euclidean norm in $\mathbb{R}^2$ (actually, any norm will do). Both the continuity and the bound come from the simple observation that on any $i$-characteristic, $R_i$ is constant and $R_j$ is a Lipschitz continuous function of the tangent inclination (with Lipschitz constant depending solely on the system $\Theta$).

For our next step we consider translates $U_1' = U_1 - m_1$ and $U_2' = U_2 - m_2$ which have the common point $0$. Note that $U_1'$ and $U_2'$ are orthogonal to each other at $0$ and so can be used as initial curves for a characteristic initial value problem, $U_k'$ being a $k$-arc. For $\epsilon$ sufficiently small (so that $U_1$ and $U_2$ are virtually orthogonal to the horizontal line $\Im\{z\} = -\epsilon$ at $e_1$ and $e_2$, respectively) convexity considerations easily show that for any admissible corner value $\rho = R(0)$ the mapping $\zeta(1, U_1', U_2', \rho, t)$ is one-to-one on the corresponding characteristic coordinate rectangle. To see that this is indeed the case, note that since the change in $\theta$ from the left end of $U_1'$ to the right end of $U_2'$ is $-(\frac{3\pi}{2} + 2\epsilon)$ (that is, the change in $\theta$ along $U_1'$ + the change in $\theta$ along $U_2'$, both from left to right, is $-(\frac{3\pi}{2} + 2\epsilon)$), the total change along the top two sides of the characteristic quadrilateral is also $-(\frac{3\pi}{2} + 2\epsilon)$.) The resulting solution $R_\rho(z) = R(1, U_1', U_2', \rho, z)$ and the quadrilateral $Q_\rho = Q(1, U_1', U_2', \rho)$ itself depend solely on $\rho = R(0)$. Note also that because $U_1', U_2' \in \mathcal{CL}$, the other two sides of $Q_\rho$ also belong to $\mathcal{CL}$.

Next we next show how we can produce curves $X_k \in \mathcal{S}(\epsilon, t_k, \alpha_k)$, such that $B(X_k \cap \mathbb{H}) \supset Q_\rho \cap \mathbb{H}$, and a solution in $B(X_k \cap \mathbb{H}) \setminus (Q_\rho \cap \mathbb{H})$ which gives a $C^\infty$ extension of $R_\rho$ for which $X_k$ is a one-sided $k$-arc which satisfies the OSC. We show how to get $X_2$, the case of $X_1$ being identical apart from obvious details.

First we extend $U_2'$ on the right by continuously adjoining a circle of radius $\frac{\epsilon}{10}$ lying in the open lower half-plane by using $F(u_2 - m_2, \lambda(U_2) - \delta, \delta, \frac{10}{\epsilon}, -(\frac{\pi}{2} + \epsilon))$, where $u_2$ parametrizes $U_2$ from left to right (i.e., with $u_2(0) = m_2$) and where $\delta$ is so small that $u_2([\lambda(U_2) - \delta, \lambda(U_2)])$ is a line segment. This in turn allows us to extend the other 2-side of $Q_\rho$ rightwards and eventually downwards (in light of the quasi-HP property) until it gets to a point $p$ for which $\Im\{p\} = -\frac{\epsilon}{2}$; call this extended 2-side $T$. The corresponding solution is defined in the simply covered domain $E = E(U_1' \cup U_2', T')$ for the subarc $T'$ of $T$ which joins a point of $\mathbb{R}$ to the left of $U_1'$ to a point of $\mathbb{R}$ to the right of $U_2'$. This means that if we change $T$ below $\mathbb{R}$ to form a new curve $T_2$, which we subsequently use as the 2-arc for a characteristic initial value problem with a straight initial 1-arc outside of $E$, then the solution is compatible with the restriction to $E$ of the solution we have so far. Using an appropriate instance of the operator $F$, we obtain an extension $T_2$ of $T$ by smoothly adjoining to $T$ at its right end a segment the argument of whose tangent angle is $\leq -\pi/2 - \epsilon$. Obviously $T_2 \in \mathcal{CL}$. It is also clear from the construction that $T_2$ depends solely on $\rho = R(0)$ and that all derivatives of its arc length parametrization depend continuously on $\rho$. This outside curve $T_2$ clearly satisfies the hypotheses (5.4) and (5.5) of the straight line characteristic initial value construction that produces an almost semicircular parallel characteristic. Let $S$ be the resulting "almost semicircle", the points at which it meets $\mathbb{R}$ being $q_1 < q_2$. By choosing $r_0$ sufficiently large we can be certain that the interior angles of this $S$ at $q_1$ and $q_2$ are within $\frac{\epsilon}{2}$ of $\frac{\pi}{2}$, and simple compactness considerations show that there is a single value of $r_0$ that is sufficiently large to have this property for all admissible corner values $\rho = R(0)$. Let $p_1$, $p_2 \in T_2$ lie on the straight 1-characteristics terminating at $q_1$ and $q_2$, respectively. It is clear that for the solution our process has generated $R(p_1)$ and $R(p_2)$ are continuous functions of $\rho$ which, moreover, satisfy bounds of the form (5.8). The solution is defined and $C^\infty$ in the part of $\mathbb{H}$ between $U_1' \cup U_2'$ and $S$, and $S$ is a 2-characteristic arc satisfying the OSC. In addition, compactness considerations show that there are absolute constants $B_1, B_2$ such that $S$ lies in $N(0, B_2) \setminus \overline{N(0, B_1)}$. Finally, we smoothly add segments to $S$ to obtain $S'$, so that for appropriate $\xi = \xi(\rho)$, $X_2 = \xi S' \in \mathcal{S}(\epsilon)$. If $x_l$ is the left endpoint of $X_2$, then $g_2(\rho) = R(x_l)$ is a continuous function of $\rho$, and similarly for the value $g_1(\rho)$ of $R$ at the left endpoint of the analogously constructed $X_1$. Here again one easily sees that in the case of $X_2$

(5.9) $\qquad \theta(g_2(\rho)) - \theta(\rho) = \theta(x_l) - \theta(0) = \frac{3\pi}{4} + \epsilon - \frac{\pi}{2} = \frac{\pi}{4} + \epsilon$

and in the case of $X_1$

(5.10) $\quad \theta(g_1(\rho)) - \theta(\rho) = \theta(x_l) - \theta(0) = -(\frac{3\pi}{4} + \epsilon) + (\frac{3\pi}{2} + 2\epsilon) = \frac{3\pi}{4} + \epsilon.$

Let $U = U_1' \cup U_2'$. By the universal property of $U_1'$, for appropriate $t^{(2)}$ and $\alpha^{(2)}$ we can use $t^{(2)} X_2 + \alpha^{(2)}$ as the 2-arc under $U_1'$ (see property (*vii*) in the discussion of the universal arcs given above), and similarly we can use some $t^{(1)} X_1 + \alpha^{(1)}$ as the 1-arc under $U_2$. These $t^{(k)} X_k + \alpha^{(k)}$ in turn come from $t^{(k)} U + \alpha^{(k)}$. If we perform the construction of $X_k$ using the value $\rho^{(k)}$ of $R$ at the center $\alpha^{(k)}$ of $t^{(k)} U + \alpha^{(k)}$, then the value of $R$ at the left endpoint of $X_k$ is $g_k(\rho^{(k)})$, so that the value of the solution $R$ (corresponding to (*vii*) of the properties of the universal arcs) in the part of the $B(U'_j)$ between $t^{(i)}U + \alpha^{(i)}$ and $U'_j$ at 0 is $f_j(g_i(\rho^{(i)}))$. We claim that there is a $\rho^{(i)}$ for which $f_j(g_i(\rho^{(i)}))$ is any given admissible corner $\rho_0$ for the characteristic initial value problem with 1- and 2-arcs $U_1$ and $U_2$. To see that such a $\rho^{(i)}$ exists in the case that $i = 2$, $j = 1$, for example, note that from (5.6) and (5.9) we have that

$$\theta(f_1(g_2(\rho^{(2)}))) - \theta(\rho^{(2)})$$
$$= \theta(f_1(g_2(\rho^{(2)}))) - \theta(g_2(\rho^{(2)})) + \theta(g_2(\rho^{(2)})) - \theta(\rho^{(2)})$$
$$= -(\tfrac{5\pi}{4} + \epsilon) + \tfrac{\pi}{4} + \epsilon = -\pi.$$

Thus, given $\rho_0$, $\rho^{(2)}$ must lie on the curve

(5.11) $$\theta(R) = \theta(\rho_0) + \pi.$$

Now, the level curves $\theta(R) = \theta_0$ in the $R$-plane are the graphs $R_i = h_{i,\theta_0}(R_j)$ of monotone functions, for which $h'_{i,\theta_0}$ are uniformly bounded and uniformly bounded away from 0. Thus as $\rho^{(2)}$ moves along the curve (5.11) it follows from the continuity of $f_1$ and $g_2$ and the bound (5.8) and the corresponding bound for the $g_k$ that $f_1(g_2(\rho^{(2)}))$ traces out the entire curve $\theta(R) = \theta(\rho_0)$, so that there is indeed a (unique) $\rho^{(1,2)}$ for which $f_1(g_2(\rho^{(2)})) = \rho_0$. The case of $f_2(g_1(\rho^{(1)})) = \rho_0$ is handled in the same manner, using (5.7) and (5.10) instead of (5.6) and (5.9).

It is now clear that we can iterate this construction to produce a solution $R^* = (R_1^*, R_2^*)$ in $B(U'_1 \cap \mathbb{H}) \cup B(U'_2 \cap \mathbb{H})$, and in fact a solution in all of $\mathbb{H}$ by appropriately extending the solution we have in $X_1$ (or in $X_2$, for that matter). More specifically, the solution in $B(U'_1 \cap \mathbb{H}) \cup B(U'_2 \cap \mathbb{H})$ is obtained by placing a suitable similar copy of the form $tU + \alpha$ under each of $U'_k$ in the way described and then under each of the corresponding $tU'_k + \alpha$ smaller similar copies of the form $tU + \alpha$ (with a smaller value of $t$, of course). At the $n^{\text{th}}$ stage (where we regard the initial $U$ as corresponding to the $0^{\text{th}}$ stage) we have $2^n$ disjoint similar copies $U_l^{(n)}$, $1 \leq l \leq 2^n$ of $U$. Let $I_l^{(n)}$ be the interval of $\mathbb{R}$ joining the left and right endpoints of $U_l^{(n)} \cap \mathbb{H}$. It is an immediate consequence of the shape of the $U_k$ that (once $\epsilon$ has been fixed) there is a constant $\delta_1 = \delta_1(\Theta) > 0$ such that for all points $p \in I_l^{(n)}$

$$\lambda(\theta(R^*(\{z : \delta_1 < \arg\{z - p\} < \pi - \delta_1\}, z \in U_l^{(n)} \cap \mathbb{H}))) \geq \delta_1.$$

Since on $U_i$, $\theta$ and $R_j$ are bi-Lipschitz functions of each other, it follows that there are $r_1 = r_1(\Theta)$, $\delta_2 = \delta_2(\Theta)$ such that for all $p \in I_l^{(n)}$ the range of $R_k^*$ in $N(r_1^n, p) \cap \{z : \delta_1 < \arg\{z - p\} < \pi - \delta_1\}$ is an interval of length at least $\delta_2$, $k = 1, 2$. Let $M^{(n)} = \cup \{I_l^{(n)} : 1 \leq l \leq 2^n\}$. Then $M = \cap \{M^{(n)} : n \geq 1\}$ consists entirely of boundary singularities of the solution $R^*$ that we have constructed. To those who have worked with Cantor sets it is probably obvious that $M$ has positive Hausdorff dimension, but for completeness we include an appropriately modified version of the argument given by Falconer [F] for the classical excluded middle third case.

First of all, it is clear from the self-similar nature of the construction that there is a number $\gamma \in (0, 1)$ such that the minimum distance between the $2^n$ intervals making up the set $M^{(n)}$ is bounded below by $\gamma^n$. We show that the $\tau$-dimensional Hausdorff-measure of $M$ is positive, where $\tau$ is defined by

$$\gamma^\tau = \tfrac{1}{2}.$$

Let $\{G_k\}$ be an open cover of $M$, which we can assume to be finite since $M$ is compact. Let $\max\{\operatorname{diam}(G_k)\} < 1$. For each $k$ there is an $l = l(k)$ such that
$$\gamma^{l+1} \leq \operatorname{diam}(G_k) < \gamma^l.$$
From the definition of $\gamma$ it follows that such a $G_k$ can have a nonempty intersection with at most one of the intervals that make up $M^{(l)}$. Consequently, for $p \geq l(k)$ at most $2^{p-l(k)}$ of the intervals making up $M^{(p)}$ can have a nonempty intersection with such a $G_k$. For each $p \geq l = l(k)$ we therefore have
$$2^{p-l} = \frac{2^p}{2^l} = \frac{2^p(\gamma^{l+1})^\tau}{2^l \gamma^{\tau l} \gamma^\tau} \leq \frac{2^p(\operatorname{diam}(G_k))^\tau}{(2\gamma^\tau)^l \gamma^\tau}.$$
Let $p$ be so large that $\min\{\operatorname{diam}(G_k)\} \geq \gamma^{p+1}$, so that $p \geq l(k)$ for all $k$. Since there are $2^p$ intervals in $M^{(p)}$, if we denote by $N_k$ the number of intervals touched by $G_k$, then
$$2^p \leq \sum_k N_k \leq \sum_k \frac{2^p(\operatorname{diam}(G_k))^\tau}{(2\gamma^\tau)^{l(k)} \gamma^\tau} = \sum_k \frac{2^p(\operatorname{diam}(G_k))^\tau}{\gamma^\tau},$$
so that
$$\sum_k (\operatorname{diam}(G_k))^\tau \geq \gamma^\tau = \tfrac{1}{2}.$$
This means that the $\tau$-dimensional Hausdorff measure of $M$ is at least $\tfrac{1}{2}$, so that $M$ has Hausdorff dimension at least $\tau$.

## 6. A Few Concluding Remarks

We briefly discuss some of the issues and problems suggested by the foregoing. In the first place, it would be interesting to determine whether Corollary 4.2 is true for genuinely nonlinear systems (see the definitions of the terms "system" and "genuinely nonlinear" between relations (1.4) and (1.5)) which are not necessarily normal. Probably, though, some hypothesis in the spirit of (*ii*) of Definition 1.1 as well as some bound like
$$\epsilon < |\theta_1(R) - \theta_2(R)| < \pi - \epsilon, \text{ for all } R \in \mathbb{R}^2$$
is necessary so that, using the approach of Section 5, or otherwise, one might try to construct a solution of a $2 \times 2$ genuinely nonlinear hyperbolic system not satisfying one or the other or both of these conditions and which has a set of boundary singularities of Hausdorff dimension 1. In a wider context one can ask if Corollary 4.2 has any counterparts for an appropriate class of sufficiently nonlinear $n \times n$ planar hyperbolic systems. In reference to normal systems, our analysis leaves open the question of whether in a half-plane $\mathbb{H}$, for example, there can be a solution for which the set of boundary singularities of type 1 has positive Hausdorff dimension. Corresponding to such a solution there would have to be set $A \subset \partial G$ of positive Hausdorff dimension such that for each $a \in A$ there is a characteristic $C_a$ exiting at $a$ but for which hypothesis (3.22) of Proposition 3.31 does not hold. Also open is the question of whether the word "nontangential" is necessary in the conclusion of the Main Theorem 2.2.

As we shall show elsewhere, the ideas of Section 5 can be used to construct in an arbitrary Jordan domain cps-mappings with arbitrary (distinct) principal stretch factors which have infinitely many isolated singularities, and in fact these singularities can be of spiral type (see [G3] for the classification of isolated singularities of cps-mappings). This raises the question of the distribution of such singularities, and in this regard we conjecture that there is some absolute constant $\gamma > 1$ such that if $\{a_n\}$ is the sequence of isolated singularities of any cps-mapping in any smoothly bounded Jordan domain for which $\{\text{dist}(a_n, \partial G)\}$ is nonincreasing, then

(6.1) $$\sum_n \text{dist}(a_n, \partial D)\gamma^n < \infty.$$

Note that this was shown with $\gamma = 1$ in [G3, Corollary 4.1]. More generally, there are other $2 \times 2$ genuinely nonlinear systems for which there exist corresponding unambiguously defined nets of characteristics which have isolated singularities, and for any such system one could attempt to obtain a classification of such singularities along the lines of [G3] and seek a bound of type (6.1) on their density.

E-mail address: jgevirtz@bsu.edu